\documentclass[a4paper,12pt]{article}

\usepackage{color}  \newcommand{\rot}[1]{{\color{black} #1}}     
\usepackage{amsmath,latexsym,amssymb}
\usepackage[german,english]{babel}
\newtheorem{thm}{Theorem}
\newtheorem{prop}{Proposition}[section]
\newtheorem{lemma}[prop]{Lemma}

\newcommand{\proof}[1][]{{\it Proof#1: }}
\newcommand{\qed}[1][3mm]{\hspace*{\fill} $\Box$ \vspace{#1}}

\newcommand{\ZZ}{{\mathbb Z}}
\newcommand{\CC}{{\mathbb C}}
\newcommand{\RR}{{\mathbb R}}
\newcommand{\HH}{{\mathbb H}}

\newcommand{\PP}{{\mathbf P}}
\newcommand{\FF}{{\mathbf F}}

\newcommand{\del}{\partial}
\newcommand{\ON}{\operatorname}

\renewcommand{\a}{\alpha}
\renewcommand{\b}{\beta}
\newcommand{\cg}{\gamma}

\renewcommand{\d}{\delta}

\newcommand{\varth}{\vartheta}

\renewcommand{\l}{\lambda}

\newcommand{\s}{\sigma}
\newcommand{\oo}{\omega}

\newcommand{\tto}{\longrightarrow}

\newcommand{\surj}{\to\!\!\!\!\!\!\!\!\tto}

\newcommand{\inv}{^{^{-1}}}

\newcommand{\afami}{{\cal A}}
\newcommand{\bfami}{{\cal B}}
\newcommand{\cfami}{{\cal C}}
\newcommand{\dfami}{{\cal D}}

\newcommand{\ffami}{{\cal F}}

\newcommand{\gfami}{{\cal G}}
\newcommand{\hfami}{{\cal H}}

\newcommand{\mfami}{{\cal M}}

\newcommand{\ofami}{{\cal O}}

\newcommand{\rfami}{{\cal R}}

\newcommand{\ufami}{{\cal U}}
\newcommand{\vfami}{{\cal V}}
\newcommand{\wfami}{{\cal W}}

\newcommand{\ibold}{{\boldsymbol{i}}}
\newcommand{\jbold}{{\boldsymbol{j}}}
\newcommand{\kbold}{{\boldsymbol{k}}}

\newcommand{\xbold}{{\boldsymbol{x}}}
\newcommand{\nubold}{{\boldsymbol{\nu}}}

\renewcommand{\xbold}{x}
\renewcommand{\nubold}{\nu}

\newcommand{\cutoff}[1]{}
\newcommand{\labell}[1]{\label{#1}}
\newcommand{\reff}[1]{\ref{#1}}

\newcommand{\pgl}{\ON{PGL}}
\newcommand{\gl}{\ON{GL}}
\newcommand{\slz}{\ON{SL_2\ZZ}}

\newcommand{\dip}{p_{n,d}}
\newcommand{\bip}{q_{n,d}}

\newcommand{\lep}{\ell_{n,d}}
\newcommand{\degz}{\deg_z}
\newcommand{\wdeg}{\ON{w-deg}}
\newcommand{\degv}{\deg_v}

\newcommand{\dnd}{\dfami_{n,d}}
\newcommand{\hnd}{\hfami_{n,d}}
\newcommand{\mnd}{\mfami_{n,d}}
\newcommand{\und}{\ufami_{n,d}}
\newcommand{\vnd}{V_{n,d}}
\newcommand{\wnd}{\wfami_{n,d}}
\newcommand{\xnd}{X_{n,d}}
\newcommand{\ynd}{Y_{n,d}}

\newcommand{\ind}{I_{n,d}}

\newcommand{\vernd}{\vfami_{n,d}}

\newcommand{\ncount}{\kappa}

\newcommand{\kord}{\prec_\ncount}

\newcommand{\vplus}{v_\Sigma} 

\newcommand{\dd}{{3d}}

\newcommand{\dms}{\textstyle}

\begin{document}

{\noindent \Large \bf Fundamental groups of moduli stacks
of smooth\\[3mm] Weierstrass fibrations\\[6mm]}

{\noindent Michael L\"onne, Universit\"at Bayreuth,\ \today\\[2mm]}


\begin{abstract}
We give finite presentations for the fundamental group of moduli stacks
of smooth Weierstrass curves over $\PP^n$ which extend
the classical result for elliptic curves to positive dimensional base.
We thus get natural generalisations of $\slz$ presented in terms of
$\left(\begin{smallmatrix} 1&1\\0&1\end{smallmatrix} \right)$,
$\left(\begin{smallmatrix} \phantom{-}1&0\\-1&1\end{smallmatrix} \right)$
and pave the way to understanding the fundamental group
of moduli stacks of elliptic surfaces in general.

Our approach exploits the natural $\ZZ_2$-action on Weierstrass curves and
the identification of $\ZZ_2$-fixed loci with smooth hypersurfaces in an
appropriate linear system on a projective line bundle over $\PP^n$.
The fundamental group of the corresponding discriminant complement
can be presented in terms of finitely many generators and relations
using methods in the Zariski tradition,
which were successfully elaborated in \cite{habil}.
\end{abstract}

%
%
%
%
%

\section{Introduction}

Our primary objects are hypersurfaces of the ruled
manifold
$\xnd=\PP\left(\ofami_{\PP^n}(d)\oplus\ofami_{\PP^n}\right)$
in the linear system $|3\s_0|$,
where $\s_0$ denotes the divisor on $\xnd$ defined
by the zero section of $\ofami_{\PP^n}(d)$.
They are assembled in the universal hypersurface $\hnd$
which is a hypersurface
in $\xnd\times\PP \vnd$, 
$\vnd=\Gamma(\xnd,\ofami_{\xnd}(3\s_0))$.

Upon a choice of homogeneous coordinates $y,y_0$ on a fibre and
$x_0,...,x_n$ on the base, $\vnd$ is identified with the polynomials
of $\CC[y_0,y,x_0,...,x_n]$ which have degree $3$ in the variables
$y_0,y$ and weighted degree $3d$ in the variables $y,x_i$ of
weights $d$ and $1$ respectively.

On $\vnd$ we introduce coordinates $u_\nu$ with respect to the
monomial basis.
Since each monomial in $\vnd$ is uniquely determined by the
exponents of the $x_i$, each coordinate $u_\nu$ is unambiguously
specified by a multiindex $\nu\in\{(\nu_1,...\nu_n)\,|\,\nu_i\geq0,
|\nu|=\sum_i\nu_i\leq3d\}$.
The equation of $\hnd$ then reads in multiindex notation, $\xbold^\nubold=x_0^{\nu_0}\cdots x_n^{\nu_n}$,
\begin{equation}
\label{subruled}
u_{\mathbf 0}y^3 
\,+\,\sum_{|\nu|=d} u_\nubold y_0y^2\xbold^\nubold
\,+\,\sum_{|\nu|=2d} u_\nubold y_0^2y\xbold^\nubold
\,+\,\sum_{|\nu|=3d} u_\nubold y_0^3\xbold^\nubold
\end{equation}

Its projection to the factor $\PP\vnd$ has singular values
precisely along the discriminant
\begin{eqnarray*}
\dnd & = &
\{ u\in\PP \vnd\:|\:\hfami_u \text{ is singular}Ê\,\}
\end{eqnarray*}
which is the union of the hyperplane $\{u_{\mathbf 0}=0\}$
and the projective
dual of weighted projective space $\PP^{n+1}_{d,1,1,...,1}$ given
as the image of $\xnd$ under the projective morphism defined
by the base point free linear system $|3\s_0|$.

The problem we want to address in the first stage is to give a
geometrically distinguished finite
presentation of the fundamental group
of the complement $\und$ of $\dnd$.

It may be viewed as a special instance
of the vastly open problem posed by Dolgachev and Libgober, \cite{dl},
to determine the fundamental group of the discriminant complement
of any (complete) linear system.

They handle the case of linear systems of elliptic curves on $\PP^2$
and $\PP^1\times\PP^1$ as well as linear systems on curves, 
but actually the first result of that kind is due to Zariski
who considered the complete linear systems on $\PP^1$, which
in degree $l=3$ may be viewed as the $n=0$ analogue of
our set up.

\begin{thm}[Zariski \cite{z} and Fadell, van Buskirk \cite{fvb}]
\hspace*{\fill}\phantom{x}\\
\labell{zariski}
The fundamental group $\pi_1(\ufami_{\PP^1,l})$
of the discriminant complement associated to the complete linear
system of degree $l$ on $\PP^1$
is finitely presented by generators
$\s_1,...,\s_{l-1}$ and relations
\begin{enumerate}
\item
$\s_i\s_j=\s_j\s_i$, if $|i-j|\geq2$, $1\leq i,j<l$,
\item
$\s_i\s_{i+1}\s_i=\s_{i+1}\s_i\s_{i+1}$, if $\,1\leq i<l-1$,
\item
$\s_1\cdots\s_{l-2}\s_{l-1}\s_{l-1}\s_{l-2}\cdots\s_1=1$.
\end{enumerate}
\end{thm}
We have previously extended these results to all complete linear
systems on projective spaces, cf.\ \cite{pdisc}, and we provide now
a series of examples of linear systems on ruled manifolds
(The claim will be made more explicit later in the paper,
cf.\ section \reff{results}.):

\begin{thm}
\label{theorem2}
There is a graph $\Gamma_{n,d}=(\vernd,E_{n,d}\subset{\vernd}^2)$ with vertex set
$\vernd$ linearly ordered by $\prec$,
such that $\pi_1(\und)$ is generated by elements
$$ 
t_\ibold,\quad \ibold\in \vernd
$$
and a complete set of relations is provided by:
\begin{enumerate}
\item
$ t_\ibold t_\jbold= t_\jbold t_\ibold$
for all $(\ibold,\jbold)\not\in E_{n,d}$,
\item
$ t_\ibold t_\jbold t_\ibold= t_\jbold t_\ibold t_\jbold$
for all $(\ibold,\jbold)\in E_{n,d}$,
\item
$ t_\ibold t_\jbold t_\kbold t_\ibold= t_\jbold t_\kbold t_\ibold t_\jbold$
for all $\ibold\prec\jbold\prec\kbold$
such that $(\ibold,\jbold),(\ibold,\kbold),(\jbold,\kbold)\in E_{n,d}$,
\item
for all $\ibold\in \vernd$\\[-2mm]
$$
 t_\ibold
 \left( t_\ibold\inv \prod_{\ibold\in\vernd^\prec} t_{\ibold}\right)^{\dd-1}
\,=\quad
\left( t_\ibold\inv \prod_{\ibold\in\vernd^\prec} t_{\ibold}\right)^{\dd-1}
t_\ibold
$$
\end{enumerate}
\end{thm}
%
In the second stage we are interested in hypersurfaces given by a
Weierstrass equation.
On the $\PP^2$-bundles $\ynd=
\PP\left(\ofami_{\PP^n}\oplus\ofami_{\PP^n}(2d)
\oplus\ofami_{\PP^n}(3d)\right)$
we introduce homogeneous fibre coordinates $y_0,y_1,y_2$
and homogeneous base coordinates $x_0,...,x_n$.

A Weierstrass fibration is then defined to be given by an equation
of the form
\begin{equation}
\label{wform}
y_0y_2^2\quad=\quad y^3 
\,+\,\sum_{|\nu|=2d} u_\nubold y_0^2y\xbold^\nubold
\,+\,\sum_{|\nu|=3d} u_\nubold y_0^3\xbold^\nubold
\end{equation}
with $d$ necessarily even, 
where the coefficients $u_\nu$ are coordinates of a vector subspace
$\vnd'$ of $\vnd$ generated by the elements $x^\nu,|\nu|\geq2d$.

The same equation also defines
the associated tautological Weierstrass hypersurface in $\ynd\times
\vnd'$. Its projection to $\vnd'$
has singular values along the discriminant
$$
\dnd' \quad = \quad
\dnd\,\cap\, \vnd'
$$
where $\vnd'$ is embedded into $\PP\vnd$ as the affine part
(w.r.t.\ the hyperplane $u_{\mathbf 0}=0$) of the
projective subspace generated by $\vnd'$ and $x^{\mathbf 0}$,
the polynomial with coefficient $u_{\mathbf 0}$.
\\
That property is easily checked on equations and reflects the
following fact.
A Weierstrass fibration is a double cover of a smooth surface
and therefore smooth if the $\ZZ_2$-fixed locus is. Its fixed part
off the hypersurface $y_2=0$ is always smooth, hence smoothness
is equivalent to smoothness of the restriction to $y_2=0$ which 
yield precisely the smoothness condition considered in the first part.
\\

The complement $\und'$ of $\dnd'$ in $\vnd'$ is the base
of a versal family of smooth Weierstrass fibrations. The moduli stack
$\mnd$ is in naturally obtained as the quotient of 
of $\vnd'$ by the group of automorphisms, which
is given as the direct product of the group of linear projective
transformations of the base and $\CC^*$ acting on the coordinates
$u_\nu$ with weight $|\nu|/d$.
We will show $\pi_1(\und)\cong\pi_1(\und')$ and derive
the topological fundamental group of $\mnd$ from
a homotopy exact sequence.

In this way we are able to generalise
the old result giving the fundamental group
of the moduli stack of elliptic curves, i.e.\ the $n=0$ case.

\begin{thm}
\label{elliptic}
The (orbifold) fundamental group $\slz$ of $\mfami_{0}$ 
is finitely presented as
$$
\langle\quad \s_1,\s_2\quad|\quad \s_1\s_2\s_1\quad=\quad \s_2\s_1\s_2,\quad
(\s_1\s_2)^6\quad=\quad 1\rangle
$$
\end{thm}

Of course the natural generalisation relies heavily on our theorem
\ref{theorem2}:

\begin{thm} 
There is a graph $\Gamma_{n,d}=(\vernd,E_{n,d})$ with 
linearly ordered vertex set 
$\vernd$ and bijective maps for $\ncount\in\{0,...,n\}$,
$$
\ibold_\ncount: \{1,...,2(\dd-1)^n\}\to \vernd
$$
such that the (orbifold) fundamental group $\pi_1(\mnd)$ 
is generated by elements
$t_\ibold,\,\ibold\in \vernd$
and a complete set of relations is given by $i)-iv)$ above
and two additional relations
\begin{enumerate}
\setcounter{enumi}{4}
\item
\[
\prod_{\ncount=0}^n\:\left(\prod_{m=1}^{2(\dd-1)^n} 
t_{\ibold_\ncount(m)}\right)^6
\quad=\quad 1\quad=\quad 
\left( \prod_{m=1}^{2(\dd-1)^n} t_{\ibold_0(m)}\right)^{d}
\]
\end{enumerate}
\end{thm}

Theorem \ref{elliptic} is recovered with $\Gamma$ the connected 
graph on $2$ vertices and the map
$\ibold_0$ enumerating its vertices.
In that case the sets $ii),iii)$ and $iv)$ of relations are void.
\\

We want to stress the fact that the relations $i)-iii)$ of our presentations have a distinctive
flavour since they stem from a different setting:
If we consider the Brieskorn-Pham polynomial in the variables $y,x_1,..,x_n$,
$$
f\quad=\quad y^3\,+\,x_1^\dd\,+\,\cdots\,+\, x_n^\dd,
$$
we are naturally led to consider a versal unfolding of the isolated hypersurface singularity it defines.
In fact the complement of the discriminant in the unfolding base was shown to have fundamental group
generated as in the theorem but with relations $i)-iii)$ only, \cite{habil}, in terms of a distinguished
Dynkin graph $\Gamma_{n,d}$ associated to the singularity, eg. the following and higher dimensional analoga:
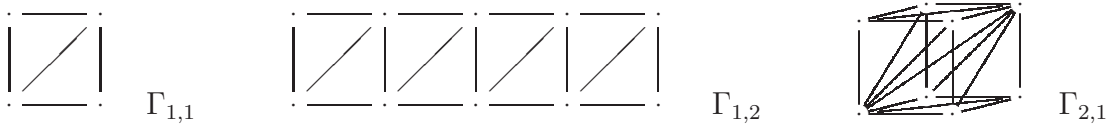
\unitlength=1.2mm

\begin{picture}(130,22)(-5,0)

\put(0,0){\begin{picture}(10,14)(5,-6)

\put(0,0){\circle*{.3}}
\put(0,10){\circle*{.3}}
\put(10,0){\circle*{.3}}
\put(10,10){\circle*{.3}}

\put(0,1.5){\line(0,1){7}}
\put(10,1.5){\line(0,1){7}}
\put(1.5,1.5){\line(1,1){7}}
\put(1.5,0){\line(1,0){7}}
\put(1.5,10){\line(1,0){7}}

\end{picture}}

\put(10,5){$\Gamma_{1,1}$}

\put(25,0){
\begin{picture}(35,14)(0,-6)

\put(0,0){\circle*{.3}}
\put(0,10){\circle*{.3}}
\put(10,0){\circle*{.3}}
\put(10,10){\circle*{.3}}
\put(20,0){\circle*{.3}}
\put(20,10){\circle*{.3}}
\put(30,0){\circle*{.3}}
\put(30,10){\circle*{.3}}
\put(40,0){\circle*{.3}}
\put(40,10){\circle*{.3}}

\put(1.5,1.5){\line(1,1){7}}
\put(11.5,1.5){\line(1,1){7}}
\put(21.5,1.5){\line(1,1){7}}
\put(31.5,1.5){\line(1,1){7}}

\put(1.5,0){\line(1,0){7}}
\put(1.5,10){\line(1,0){7}}
\put(11.5,0){\line(1,0){7}}
\put(11.5,10){\line(1,0){7}}
\put(21.5,0){\line(1,0){7}}
\put(21.5,10){\line(1,0){7}}
\put(31.5,0){\line(1,0){7}}
\put(31.5,10){\line(1,0){7}}

\put(0,1.5){\line(0,1){7}}
\put(10,1.5){\line(0,1){7}}
\put(20,1.5){\line(0,1){7}}
\put(30,1.5){\line(0,1){7}}
\put(40,1.5){\line(0,1){7}}

\end{picture}}

\put(72,5){$\Gamma_{1,2}$}

\put(87,5){\unitlength=.35pt
\begin{picture}(200,200)

\put(0,0){
\begin{picture}(100,100)
\bezier{110}(0,0)(50,0)(100,0)
\bezier{110}(0,100)(50,100)(100,100)
\bezier{110}(0,0)(0,50)(0,100)
\bezier{110}(100,0)(100,50)(100,100)
\bezier{113}(0,0)(50,50)(100,100)

\bezier{110}(0,0)(36,9)(72,18)
\bezier{110}(100,0)(136,9)(172,18)
\bezier{110}(0,100)(36,109)(72,118)
\bezier{110}(100,100)(136,109)(172,118)

\bezier{113}(20,2)(86,9)(152,16)
\bezier{113}(20,102)(86,109)(152,116)
\bezier{113}(0,0)(36,59)(72,118)
\bezier{113}(100,0)(136,59)(172,118)
\bezier{130}(0,0)(86,59)(172,118)

\put(0,0){\color{white}\circle*{20}}
\put(100,0){\color{white}\circle*{20}}
\put(0,100){\color{white}\circle*{20}}
\put(100,100){\color{white}\circle*{20}}

\put(0,0){\circle*{4}}
\put(100,0){\circle*{4}}
\put(0,100){\circle*{4}}
\put(100,100){\circle*{4}}

\end{picture}}

\put(72,18){
\begin{picture}(100,100)
\bezier{110}(0,0)(50,0)(100,0)
\bezier{110}(0,100)(50,100)(100,100)
\bezier{110}(0,0)(0,50)(0,100)
\bezier{110}(100,0)(100,50)(100,100)
\bezier{113}(0,0)(50,50)(100,100)
\end{picture}}

\put(72,18){
\begin{picture}(100,100)
\put(0,0){\color{white}\circle*{20}}
\put(100,0){\color{white}\circle*{20}}
\put(0,100){\color{white}\circle*{20}}
\put(100,100){\color{white}\circle*{20}}

\put(0,0){\circle*{4}}
\put(100,0){\circle*{4}}
\put(0,100){\circle*{4}}
\put(100,100){\circle*{4}}
\end{picture}}

\end{picture}}

\put(110,5){$\Gamma_{2,1}$}

\end{picture}
We will explain in detail in section \reff{affine} how this result can be used in the present paper.

Relations $iv)$ on the contrary are due to degenerations along the
hypersurface $x_0=0$, while those in $v)$ originate in the action
of the automorphism group. 
\\

The present paper should be viewed as a further contribution in our
ongoing project to understand families of smooth elliptic surfaces
and their monodromies, for which we have given an outline in the
introduction of \cite{bifbraid}.
In fact \cite{reg} and \cite{ireg} may be seen as a starting point, since
there we have determined the images of homological monodromy.

Moduli stacks enter the stage, since they provide the appropriate
means to study all families of a specified kind at once.
In particular all their monodromy maps should assemble into a
monodromy homomorphism defined on the topological fundamental
group of the stack.

A particular nice example -- which motivated our research --
is provided by the families of elliptic curves, where the homological
monodromies assemble into an isomorphism from the orbifold
fundamental group of the stack $\HH/\slz$
to the automorphism group $\slz$ of the first homology of a curve,
cf.\ theorem \ref{elliptic}.
Our aim is to investigate possible generalisations to the case of
families of elliptic surfaces, which we believe to be tractable and
still to exhibit many characteristic features of the surface case in
general.
\\

A major difference to the curve case is the existence of  -- at least --
three distinct moduli problems for families of elliptic surfaces
which attract our attention:
\begin{enumerate}
\item
for smooth elliptic surfaces with a section. The coarse moduli space 
has been
constructed by Miranda and Seiler as the moduli space of 
Weierstrass fibrations with at most rational double points,
cf. \cite{mi,sei1}.
\item
for smooth elliptic surfaces with a section and irreducible fibres only,
equivalently for surfaces with a smooth Weierstrass model. That case
is an instance of a moduli problem for polarised elliptic surfaces
as considered by Seiler \cite{sei2}.
\item
for smooth elliptic surfaces with a section and nodal fibres only,
which were considered in \cite{bifbraid} for the benefit of allowing
a special kind of monodromy, cf.\ below.
\end{enumerate}

To hope for as nice a result as in the elliptic curve case, we are forced
to adjust the choice of monodromy to the choice of moduli problem.
An educated guess among some natural monodromies leads to the
following tentative list:
\begin{enumerate}
\item
algebraic or geometric monodromy. It takes values in the
automorphism group of integer (co)homology respective the group
of isotopy classes of diffeomorphism.
\item
symplectic monodromy. Both the ambient space and the polarisation
may be employed to construct a symplectic connection.
The monodromy then takes values
in the group of symplectic isotopy classes of symplectomorphisms.
\item
bifurcation braid monodromy. We exploit the fact, that families
of elliptic surfaces with nodal fibres only, they naturally give rise to
continuous
families of finite sets in the base. Thus in case of regular surfaces
the monodromy takes values in the braid group of the two-sphere,
cf. \cite{bifbraid}.
\end{enumerate}
Since symplectic monodromy remains quite mysterious despite
the efforts of Seidel and others to enlighten the structure of
symplectomorphism groups we have proposed a replacement
of $ii)$ of a more topological flavour:
(More details and motivation from a comparison with symplectic
monodromy can be found in \cite{bifbraid}.)
\begin{enumerate}
\item[ii')]
braid class monodromy: Obtained from braid monodromy
by imposing just as many relations on the image of braid monodromy
as to make sure, that it is well defined on the larger moduli stack.
\end{enumerate}

In any case it is desirable
to understand the topological fundamental groups
of the moduli stacks and the target groups of the monodromies.
While our previous contributions were to monodromies in case $i)$
and $iii)$, the present paper yields the fundamental group in case
$ii)$. 

Our results also prepare the ground to handle the fundamental
group in the other cases.
To address $iii)$ we have to discard some parts of the moduli
stack. On the level of discriminant complements this corresponds
to taking the bifurcation divisor into account, the set of parameters
$u$, such that the projection of the corresponding
hypersurface $\hfami_u$ to $\PP^1$ is non-generic.

For case $i)$ on the other hand, we need to glue in
some orbifold divisor to account for
some families which are allowed in addition.
The associate coarse space is naturally the
coarse moduli space of elliptic surfaces
with a section constructed as a moduli space of Weierstrass
fibrations with at most rational double point singularities.
To construct the appropriate stack structure over that space, to get
the actual moduli stack for families of smooth elliptic surfaces,
will be the task of a forthcoming paper.
\\

Of course we can initiate an analogous program in higher
dimension. For example our new results have no dimension restriction.
Nevertheless we should note a number of potential obstacles:
\begin{enumerate}
\item
In higher dimension a generalised bifurcation monodromy can
be assigned as long as we admit only family of Weierstrass
fibrations with generic bifurcation set of their fibrations.
However this monodromy maps but to a group
detecting the braiding in $\PP^n$ of the critical loci,
which then are positive dimensional and singular, cf. the interpretation
of $\pi_1(\ufami_{\PP^n\!\!,d})$ as group of braiding in $\PP^n$ 
(\cite{pdisc}).
\item
Admitting also families of smooth Weierstrass fibrations, the need
of a bifurcation class monodromy has to be checked and if necessary
relations have to be imposed on the image of bifurcation monodromy.
\item
A suitable relation of smooth elliptic fibrations with section to
Weierstrass fibrations with mild singularities is needed for
any progress on geometric monodromy.
\end{enumerate}

\section{the discriminant polynomial}

\newcommand{\lsys}{|\ofami_{\xnd}(3\s_0)|}

The aim of this section is to gain a better insight into the geometric
properties of $\dnd$ and $\und$.
Upon identification of the complement of the hyperplane
$u_{\mathbf 0}=0$ in $\PP \vnd$ with the affine hyperplane
$u_{\mathbf 0}=1$ in $\vnd$, $\und$ is the complement of a
hypersurface given by a polynomial 
$$
\dip \in
\CC\big[ u_\nu\,\big|\, |\nu|\in\{d,2d,3d\}\big].
$$
We distinguish the parameter $z=u_{\dd,0,..,0}$ (and emphasize the
distinction by the new notation from now on). With respect to
the parameter $z$
we define the discriminant polynomial $\bip=discr_z\dip$ of
$\dip$ and its leading coefficient $\lep$, which together with $\dip$
will be the targets of our ensuing investigations.

For convenience we recall some topological Euler numbers:
$$
\begin{array}{rcll}
e_n &=& n+1, &
\text{of complex projective space }\PP^n,\\[3mm]
e_{n;d} &=&  n+1+\frac{\dms(1-d)^{n+1}-1}{\dms d}, &
\text{of smooth hypersurfaces of degree $d$ in }
\PP^n,\\[4mm]
e_{n;d,d}  &=& 
n+1 + (n-1)(1-d)^n &
\text{of smooth complete intersections of two}\\[1mm]
&&
\qquad+2\frac{\dms(1-d)^{n}-1}{\dms d}, &
\quad\text{hypersurfaces of degree $d$ in }\PP^n,\\[4mm]
e_{_H} &=&  3 e_n - 2 e_{n;\dd}, &
\text{of a smooth $H$ in }
|\ofami_{\xnd}(3\s_0)|,\\[3mm]
e_{_{H\cap H'}} &=& 
3 e_{n;\dd} - 2 e_{n;\dd,\dd}, &
\text{of a smooth intersection of two}\\[1mm]
&&& \quad \text{hypersurfaces in }
|\ofami_{\xnd}(3\s_0)|.
\end{array}
$$
The last two formulas are immediate from the fact that
there is a smooth hypersurface $H$
which is a triple cover over $\PP^n$ totally branched over
a hypersurface of degree $\dd$, resp.\ an intersection of two such
hypersurfaces which is a triple cover over a hypersurface
of degree $\dd$ totally branched over the intersection with another
such hypersurface in $\PP^n$.

\begin{lemma}
\labell{degp}
The discriminant polynomial $\dip$ is of degree
$$
\deg \dip\quad = \quad 2(n+1)(\dd-1)^n.
$$
\end{lemma}

\proof
The degree of $\dip$ is by definition the number intersection points
of its zero set with a generic affine line, hence the number of
singular hypersurfaces of the corresponding affine pencil.

The hypersurfaces of such a pencil are contained in an open part of
$\xnd$ isomorphic to $\ofami_{\PP^n}(d)$.
We obtain the degree of $\dip$ computing the Euler number
$n+1=e_\CC e_n$ of $\ofami_{\PP^n}(d)$ from a decomposition
into constructible strata with respect to a generic affine pencil:

First
the set of points, which belong to fibres intersecting the base locus.
It is a $\CC$-fibre space over a degree $\dd$ hypersurface of $\PP^n$
and has Euler number $e_{n;\dd}$.

Second the set of points not in the first set, which belong to singular 
hypersurfaces.
This set consists of $\deg\dip$ hypersurfaces, each regular except
for a single ordinary double point and deprived of the base
locus of the pencil. Its Euler number is therefore
$\deg\dip(3 e_n-2e_{n;\dd}-(-1)^{n}-(3 e_{n;\dd}-2 e_{n;\dd,\dd}))$.

Third the set of points not in the first set, which belong to smooth
hypersurfaces.
It consists of a smooth family of smooth hypersurfaces
each deprived of the base locus of the pencil over the affine
line punctured at the $\deg\dip$ parameters of 
singular hypersurfaces. The Euler number is therefore
$(1-\deg\dip)(3 e_n-2e_{n;\dd}-(3 e_{n;\dd}-2 e_{n;\dd,\dd}))$.

If we equate the sum of their Euler numbers with $n+1$ and use
the numerical values provided above, we get by a straightforward
calculation
\begin{eqnarray*}
n+1 & = & e_{n;\dd}-\deg\dip(-1)^n
+e_\CC(3 e_n-2e_{n;\dd}-(3 e_{n;\dd}-2 e_{n;\dd,\dd}))\\
\iff \deg\dip & = & (-1)^n(2e_n-4 e_{n;\dd}+2 e_{n;\dd,\dd})\\
& = & 2(n+1)(\dd-1)^n\\[-14mm]
\end{eqnarray*}
\qed

\begin{lemma}
\labell{degzp}
The discriminant polynomial $\dip$ as a polynomial in the
coefficient $z$ of $x_0^\dd$ only has degree
$$
\deg_z\dip\quad=\quad2(\dd-1)^n.
$$
\end{lemma}

\proof
The degree $\deg_z\dip$ is equal to the number of singular hypersurfaces
in a generic affine pencil with varying part $z x_0^\dd$.

We obtain the degree computing the Euler number
$n+1=e_\CC e_n$ of $\ofami_{\PP^n}(d)$ from a decomposition
into three constructible strata with respect to that pencil:

First the set of points, which belong to fibres intersecting the base locus.
It is a $\CC$-fibre space over the hyperplane $x_0=0$ of $\PP^n$
and has Euler number $e_{n-1}$.

Second the set of points not in the first set, which belong to singular 
hypersurfaces.
This set consists of $\deg_z\dip$ hypersurfaces, each regular except
for a single ordinary double point and deprived of the base
locus of the pencil. Its Euler number is therefore
$\deg_z\dip(3 e_n-2e_{n;\dd}-(-1)^{n}-(3 e_{n-1}-2 e_{n-1;\dd}))$.

Third the set of points not in the first set, which belong to smooth
hypersurfaces.
It consists of a smooth family of smooth hypersurfaces
each deprived of the base locus of the pencil over the affine
line punctured at the $\deg_z\dip$ parameters of 
singular hypersurfaces. The Euler number is thus
$(1-\deg_z\dip)(3 e_n-2e_{n;\dd}-(3 e_{n-1}-2 e_{n-1;\dd}))$.

If we equate the sum of their Euler numbers with $n+1$ and use
the numerical values provided above, we get by a straightforward
calculation
\begin{eqnarray*}
n+1 & = & e_{n-1}-\deg_z\dip(-1)^n
+e_\CC(3 e_n-2e_{n;\dd}-(3 e_{n-1}-2 e_{n-1;\dd}))\\
\iff \deg_z\dip & = & (-1)^n(2e_n-2 e_{n-1}-2 e_{n;\dd}+2 e_{n-1;\dd})\\
& = & 2(\dd-1)^n\\[-14mm]
\end{eqnarray*}
\qed[8mm]

To cope with their r\^ole in the following discussion we 
introduce the shorthand $u_\nu'$ for the parameters
of monomials $x^\nu$ not containing $x_0$.

\begin{lemma}
\labell{dirr}
The discriminant $\dnd$ defined by the polynomial $\dip$ is irreducible.
\end{lemma}

\proof
In case $n=0$ the discriminant is the cuspidal cubic, hence irreducible.
In case $n>0$ the fibres of $\dnd$ under projection to the
variables $u_\nu'$ are generically irreducible.

If $\dfami$ were reducible one component thus had to be a union of
fibres and therefore had to be defined by a polynomial $g$ in
$\CC[u_\nu']$.
Of course $g$ must be a factor of each coefficient of $\dip$ considered
as a polynomial in $z$, in particular of the leading coefficient $\lep$.
By the induction hypothesis we conclude that $g$ is a factor of
$p_{n-1,d}$. But since there are parameter points on the zero set
of $p_{n-1,d}$, which correspond to smooth hypersurfaces,
it does not belong to $\dnd$. Hence there is no fibral component
of the discriminant, so we conclude that the discriminant
coincides with its irreducible vertical component.
\qed

\begin{lemma}
\labell{leadc}
The discriminant polynomial $\dip$ as a polynomial in the coefficient
$z$ of $x_0^d$ only has leading coefficient
$$
\lep \quad=\quad p_{n-1,d}^{\dd-1}.
$$
\end{lemma}

\proof
For any pencil of hypersurfaces defined by polynomials with varying
part $zx_0^\dd$
we can compute the degree $\deg_z$ of the discriminant
polynomial in $z$ by an evaluation of topological Euler numbers again.

We consider the decomposition of $\ofami_{\PP^n}(d)$
into three constructible strata with respect to the given pencil:

First the set of points, which belong to fibres intersecting the base locus.
It is the line bundle over the hyperplane $x_0=0$ of $\PP^n$
and has Euler number $e_{n-1}$.

Second the set of points not in the first set, which belong to singular 
hypersurfaces.
This set consists of $\deg_z$ hypersurfaces, each regular except
for a single ordinary double point and deprived of the base
locus $Bs$ of the pencil. Its Euler number is therefore
$\deg_z(3 e_n-2e_{n;\dd}-(-1)^{n}-e_{Bs})$.

Third the set of points not in the first set, which belong to smooth
hypersurfaces.
It consists of a smooth family of smooth hypersurfaces
each deprived of the base locus $Bs$ of the pencil over the affine
line punctured at the finitely many parameters of 
singular hypersurfaces. The Euler number is thus
$(1-\deg_z)(3 e_n-2e_{n;\dd}-e_{Bs})$.

If we equate the sum of their Euler numbers with $n+1$ and use
the numerical values provided above, we get by a straightforward
calculation
\begin{eqnarray*}
n+1 & = & e_{n-1}-\deg_z(-1)^n
+e_\CC(3 e_n-2e_{n;\dd}-e_{Bs})\\
\iff \deg_z & = & (-1)^n(2e_n-2 e_{n-1}-2 e_{Bs})
\end{eqnarray*}
where the base locus $Bs$ and $\deg_z$ depend on the pencil.

Hence the degree $\deg_z$ drops if and only if the Euler number of
$Bs$ differs by a positive multiple of $(-1)^{n-1}$ as compared to 
the Euler number of the base locus for a generic pencil.
That change occurs if and only if $Bs$ is singular, which is a condition on the restriction to the hyperplane $x_0=0$ of $H_0$.

We conclude that the degree $\deg_z$ drops if and only if the discriminant
polynomial in the appropriate variables -- the parameters $u_\nu'$
of monomials not containing $x_0$ -- vanishes.
Therefore the reduced equation for the zero set of $\lep$ is $p_{n-1,d}$, which is irreducible by the irreducibility of the discriminant.

We find the multiplicity of $p_{n-1,d}$ in $\lep$ by a comparison of degrees:
By lemma \ref{degp} the polynomial $\dip$ is homogeneous of 
degree $2(n+1)(\dd-1)^n$ and it has to match the sum of the homogenous degree
of $\lep$, which is a multiple of $\deg p_{n-1,d}=2n(\dd-1)^{n-1}$,
and $\deg_z\dip$ which is $2(\dd-1)^n$ by lemma \ref{degzp}.
So we infer our claim.
\qed

\begin{lemma}
\labell{copc}
The discriminant polynomial $\dip$ as a polynomial in the coefficient
$z$ of $x_0^\dd$ only has coprime coefficients.
\end{lemma}

\proof
By the preceding lemma
the leading coefficient has a unique irreducible factor $p_{n-1,d}$.
So the coefficients are not coprime only if $p_{n-1,d}$ is a 
factor of each.

In that case the zero set of $p_{n-1,d}$ belongs to the discriminant and must be
equal to the discriminant since both are irreducible.

This is not true because there are singular hypersurfaces which are regular
when restricted to $x_0=0$ so contrary to our assumption the coefficients are
coprime.
\qed

\begin{lemma}
\labell{degq}
The bifurcation polynomial $\bip$ is homogeneous of degree
$$
(2n+1)2(\dd-1)^n\left(2(\dd-1)^n-1\right).
$$
\end{lemma}

\proof
The polynomial $\bip$ is obtained as the discriminant of $\dip$ with
respect to the variable $z$ and so is homogeneous itself. In fact it can be computed from the
Sylvester matrix of $\dip$ and $\del_z\dip$; up to a factor consisting
of the leading coefficient polynomial $\lep$ it is the determinant of that
matrix.

Hence it is sufficient to add the degrees along the diagonal for any reordering of the matrix.
In fact we can arrange on this diagonal $\deg_z\dip$ times the leading coefficient of degree
$\deg\dip-\degz\dip$ and $\degz\dip-1$ times the constant coefficient of degree $\deg\dip$.
By the above we have to subtract the degree of the leading coefficient and thus we get
\begin{eqnarray*}
\deg\bip & = & \degz\dip(\deg\dip-\degz\dip)+(\degz\dip-1)\deg\dip\\
&&
-(\deg\dip-\degz\dip)\\
&=&
(\degz\dip-1)(\deg\dip-\degz\dip)+(\degz\dip-1)\deg\dip\\
&=&
(\degz\dip-1)(2\deg\dip-\degz\dip)\\
&=& \left(2(\dd-1)^n-1\right)\left(4(n+1)(\dd-1)^n-2(\dd-1)^n\right)\\
&=&
(2n+1)\left(2(\dd-1)^n-1\right)2(\dd-1)^n\\[-14mm]
\end{eqnarray*}
\qed[8mm]

In the next instances we consider $\dip$ to be weighted homogeneous. The weight we assign to
each $u_\nu$ is equal to the exponent of $x_0$ in the monomial of which it is the parameter:
$$
wt(u_\nu)\quad = \quad \nu_0.
$$
In particular the weight is zero if $x_0$ does not occur, ie. in case
of parameters $u_\nu'$, and it is $\dd$ precisely 
in the case of the parameter $z$.

\begin{lemma}
\labell{wdegp}
$\dip$ is weighted homogeneous of degree
$$
\wdeg\dip\quad=\quad 2\cdot\dd(\dd-1)^n.
$$
\end{lemma}

\proof
The leading term of $\dip$ is of degree $2(\dd-1)^n$ in $z$ which is of weight $\dd$
with coefficient $\lep$, a polynomial in the $u_\nu'$. Hence $\lep$ is of
weight zero and the claim follows.
\qed

\begin{lemma}
\labell{wdegq}
$\bip$ is weighted homogeneous of degree
$$
\wdeg\bip\quad=\quad 2\cdot\dd(\dd-1)^n\left(2(\dd-1)^n-1\right).
$$
\end{lemma}

\proof
The leading coefficient of $\dip$ with respect to $z$ is of weighted degree zero,
hence $\bip$ is the determinant of $\dip$ with respect to $z$ up to a factor of
vanishing weighted degree and an argument as in lemma \ref{degq} yields
$$
\wdeg\dip\quad=\quad\wdeg\dip(\degz\dip-1).
$$
With the numerical values given in lemma \ref{wdegp} and lemma \ref{degzp} we get the claim.
\qed

We next investigate the properties of $\dip$ with respect to
the parameters $u_{\nu}$ of the monomials
$x_0^{\dd-1}x_\ncount$ only. They will be called linear coefficients
and denoted by $v_\ncount$, whenever we want to emphasize
their distinguished r\^ole.

\begin{lemma}
\labell{degvq}
$\bip$ as a polynomial in the linear coefficients $v_\ncount$ is of degree
$$
\degv\bip\quad=\quad  2\cdot\dd(\dd-1)^{n-1}\left(2(\dd-1)^n-1\right).
$$
\end{lemma}

\proof
The linear coefficients are of weight $\dd-1$, so we get the upper 
bound for $\degv\bip$ to be $\wdeg\bip$ divided by $(\dd-1)$.

The existence of at least one non-trivial coefficient
can be deduced from the special family
$$
y^3-3yx_0^{2d}+\sum a_\ncount x_\ncount ^\dd
+\sum v_\ncount x_\ncount x_0^{\dd-1}+zx_0^\dd.
$$
For all $a_\ncount =1$ and all $v_\ncount $ positive real of sufficiently distinct magnitude
$$
0<v_1\ll v_2\ll \cdots\ll v_n<1
$$
the discriminant polynomial has simple roots only
so the bifurcation polynomial for the family is non-zero.

On the other hand the bifurcation polynomial of our special family is 
weighted homogeneous again with
$\wdeg=6d(\dd-1)^n(2(\dd-1)^n-1)$ to which in fact
only the $v_\ncount $ contribute since all $a_\ncount$
have weight $0$. So from the non-triviality
above we conclude our claim.
\qed

\begin{lemma}
\labell{degcv}
Consider $\bip$ as a polynomial in the parameters $v_\ncount $ with coefficients in $\CC\big[u_\nu\big|\nu_0<\dd-1\big]$.
Given a monomial of degree $\degv\bip$ in the parameters $v_\ncount $, its coefficient $c$
in $\bip$ is a polynomial in the parameters $u_\nu'$ of monomials $x^\nu$
not containing $x_0$ and it is either zero or has degree
$$
\deg c \quad = \quad 
\left(2n(\dd-1)-1\right)2(\dd-1)^{n-1}\left(2(\dd-1)^n-1\right).
$$
\end{lemma}

\proof
The degree is just the difference between $\deg\bip$ and $\degv\bip$.
The other claim follows from the fact that the leading coefficient has weighted degree $0$, so must
be a polynomial in the weight $0$ parameters $u_\nu'$.
\qed

\begin{lemma}
\labell{coprimv}
Consider $\bip$ as a polynomial with coefficients in
$\CC\big[u_\nu\big|\nu_0<\dd-1\big]$.
Then the greatest common divisor of all these coefficient polynomials
is trivial.
\end{lemma}

\proof
By construction a polynomial $f$ belongs to the zero set of $\bip$ if and only if
at least one of the polynomials $f+zx_0^\dd, z\in\CC$ has more than a single ordinary
double point singularity
or the restriction $f|_{x_0=0}$ has more than
a single ordinary double point singularity.

If $f$ is any polynomial then some perturbation $\tilde f$
of $f$ by terms $v_\ncount x_\ncount x_0^{\dd-1}$ has the property
that $\tilde f$ has non-degenerate critical points only.
Moreover by changing the perturbation ever so slightly
we may assume, that $\tilde f$ has even no multiple critical
values. Hence $\tilde f$ belongs to the zero set of $\bip$ if
and only if the restriction $\tilde f|_{x_0=0}$ is non-generically singular.

The zero set of a common factor of all coefficients is
either a divisor or empty since $\bip$ is non-trivial.

Moreover by the preceding lemma \ref{degcv} the coefficients
of monomials in $\CC[v_\ncount ]$ of highest degree contain only
parameters $u'$ of monomials $x^\nu$ not containing $x_0$,
hence the same must be true for a common
factor of all coefficients.

Therefore if a polynomial $f$ belongs to the zero set of a common factor
then so does every perturbation $\tilde f$ as above.

Hence a polynomial can only belong to the zero set
if its restriction $f|_{x_0=0}$ is non-generically singular.
Since the set of polynomials with non-generically singular restriction
to $x_0=0$ is of codimension two, the zero set of any
common factor is empty and therefore any such common
factor is a non-zero constant.
\qed

\section{Zariski arguments}

The ideas of Zariski as elaborated by Bessis \cite{bessis} for
the affine set up provide the tool to get hold of a presentation for the
fundamental group of $\und$. 
A key notion concerns distinguished elements:

\paragraph{Definition}
An element in a fundamental group representable by a
path isotopic to the boundary of a small disc transversal to a divisor
is called a \emph{geometric element}.

A basis of the fundamental group of a punctured affine line
is called a \emph{geometric basis}, if its elements are geometric
and simultaneously representable by paths disjoint except for the
base point.
\\

We denote now by $\CC^N$ the affine parameter space complement
to the hypersurface $u_{\mathbf 0}=0$ containing divisors
$\dfami,\hat\afami,\hat\bfami$ 
given by $\dip$, $p_{n-1,d}$ and $\bip$ respectively. 
Note that for notational convenience we suppress the dependence on integers $n,d$ occasionally.

We project $\CC^N$ to $\CC^{N-1}$ along the distinguished parameter
$z$ and
get divisors $\afami,\bfami$ defined by $p_{n-1,d}$ and $\bip$ again.
By construction $\afami,\bfami$ pull back to $\hat\afami,\hat\bfami$
along the projection and $\dfami$ is finite over $\CC^{N-1}$, branches
along $\bfami$ and has $\hat\afami$ as its vertical asymptotes.

\begin{lemma}
\labell{sabd} 
Suppose $L$ is a fibre of the projection such that its intersection $\dfami_L$
with the discriminant $\dfami$ consists of $\deg_z\dip$ points,
then there is a split exact sequence
$$
1\to\pi_1(L-\dfami_L)\to\pi_1(\CC^N-\hat\afami-\hat\bfami-\dfami)
\to
\pi_1(\CC^{N-1}-\afami-\bfami)\to1.
$$
with a splitting map which takes geometric elements associated to $\bfami$
to geometric elements associated to $\hat\bfami$.
\end{lemma}

\proof
In fact over the complement of $\afami\cup\bfami$ the discriminant is a finite topological cover
and its complement is a locally trivial fibre bundle with fibre the affine line punctured at
$\deg_z\dip$ points.
The exact sequence is now obtained from the long exact sequence of that fibre bundle.
Exactness on the left follows from the fact that no free group of rank more than $1$ admits
a normal abelian subgroup.

Since the points on $\dfami$ vary continuously with the parameters
in $\CC^{N-1}-\afami$ so does a suitably chosen real upper bound on
their moduli. This bound defines a topological section inducing an
algebraic one.

Moreover it maps boundaries of small discs transversal to
$\bfami$ to boundaries of small discs transversal to $\hat\bfami$ and
disjoint to any other divisor.
(Note that the last claim does not hold for boundaries of arbitrarily small discs transversal to $\afami$.)
\qed

\rot{
We derive an immediate corollary on the level of presentations:

\begin{lemma}
\labell{pabd}  
Suppose there is a presentation for the fundamental group of the base
$$
\pi_1(\CC^{N-1}-\afami-\bfami)\quad\cong\quad
\langle r_\a | \rfami_q\rangle.
$$
in terms of geometric generators then there is a presentation
$$
\pi_1(\CC^n-\hat\afami-\hat\bfami-\dfami)
\quad\cong\quad
\langle t_i,\hat r_\a | \hat r_\a t_i\inv\hat r_\a\inv\phi_\a(t_i),\rfami_q\rangle.
$$
with elements $t_i$ of a free geometric basis for a generic fibre
$L-\dfami_L$
and $\phi_\a$, $\hat r_\a$ the monodromy automorphism, resp. the lift
associated to $r_\a$.
\end{lemma}
}

To get hold on $\pi_1(\CC^{N-1}-\afami-\bfami)$
we will exploit a further projection.

\begin{lemma}
\labell{projv}
\begin{sloppypar}
There is a linear combination $\vplus$ of the $v_\ncount$ such that the projection \mbox{$p_v:\CC^{N-1}\to\CC^{N-2}$}
along $\vplus$ has the following property:
\begin{quote}
There exists a divisor $\bar\cfami$ such that no component of its pull-back $\cfami$ to $\CC^{N-1}$
is a component of $\bfami$ and such that the induced map
$p_v|:\bfami\to \CC^{N-2}$ is a topological finite covering over the complement of $\bar\cfami$.
\end{quote}
\end{sloppypar}
\end{lemma}

\proof
For general $\vplus$
the set of singular values for the induced map $p_v|:\bfami\to \CC^{N-2}$ is a divisor $\bar\cfami$, since
we equip $\bfami$ with its reduced structure. Moreover in its
complement $\bfami$ must be a topological
fibration hence a topological covering.

Since a common component of $\bfami$ and the divisor
$\cfami$ can be detected algebraically as a nontrivial factor
of $\bip$ which is independent of the variable $\vplus$,
It suffices to show that for a suitable choice of $\vplus$ there is no such factor.

We decompose the polynomial algebra $\CC[u_\nu]$ according to the degree $\deg_v$
of each monomial considered as a monomial in the $v_\ncount$ only.
With respect to the $\deg_v$ decomposition we have the summand $q_{max}$ of $\bip$ of highest degree.
Its coefficients are in $\CC[u'_\nubold]$ by lemma \reff{degcv}.

Therefore $q_{max}$ defines a proper hypersurface in some trivial affine bundle over $\PP^{n-1}$.
If we replace the projective coordinates $v_1:v_2:...:v_n$ of $\PP^{n-1}$ by new ones
$v_\ncount'$ in such a way that a point in the complement
has coordinates $(0:...:0:1)$ in the factor $\PP^{n-1}$,
then $\bip$ is a polynomial of highest possible degree
in the variable $\vplus=v_n'$. The leading coefficient
is thus in $\CC[u_\nu']$.

With that choice, a non-trivial factor of $\bip$ independent of $\vplus$ may only depend on the $u_\nu'$.
But we know already from the proof of lemma \reff{coprimv}, that no
divisor defined in terms of the $u'_\nubold$ only can be a component of
$\bfami$, hence there is no common component of $\cfami$ and $\bfami$.
\qed

We suppose from now on a projection $p_v:\CC^{N-1}\to\CC^{N-2}$
as in lemma \ref{projv} has been fixed with $L'$ a generic fibre
and we denote by $\bar\afami$ the divisor in $\CC^{N-2}$ defined
by $p_{n-1,d}$.

\begin{lemma}
\labell{genab}  
Suppose there are geometric elements $r_a$ associated to $\afami$
and a geometric basis consisting of elements $r_b$ of $\pi_1(L'-\bfami_{L'})$
such that the $r_a$ generate $\pi_1(\CC^{N-1}-\afami)$, 
then the $r_a$ and $r_b$ together generate
$$
\pi_1(\CC^{N-1}-\afami-\bfami).
$$
\end{lemma}

\proof
Given the elements $r_a$
which generate $\pi_1(\CC^{N-1}-\afami)$
we may conclude, cf. \cite{bessis}, that there are geometric elements
$r_c$ associated to $\cfami$ which can be
taken in the complement of $\bfami$ such that together they generate
$$
\pi_1(\CC^{N-1}-\afami-\cfami)\quad\cong\quad\pi_1(\CC^{N-2}-\bar\afami-\bar\cfami).
$$
In fact the situation is similar to that of lemma \reff{sabd} since
$\CC^{N-1}-\afami-\bfami-\cfami$
fibres locally trivial over $\CC^{N-2}-\bar\afami-\bar\cfami$ with fibre
the appropriately punctured complex line $L'-\bfami_{L'}$.
In the associated short exact sequence
$$
\pi_1(L_b-\bfami_L)\tto\pi_1(\CC^{N-1}-\afami-\bfami-\cfami)\tto\pi_1(\CC^{N-2}-\bar\afami-\bar\cfami)
$$
the elements $r_b$ generate the group on the left hand side and the images
of the $r_a$ and $r_c$ generate the group on the right hand side, so together they
generate the group in the middle.
But then we may deduce, that together they generate
$$
\pi_1(\CC^{N-1}-\afami-\bfami)
$$
and that the elements $r_c$ are trivial in that group, so our claim holds.
\qed

\begin{lemma}
\labell{pad}  
Suppose there is a presentation for the fundamental group of $\CC^{N-1}-\afami$
$$
\pi_1(\CC^{N-1}-\afami)\quad\cong\quad
\langle r_a | \rfami_a\rangle.
$$
in terms of geometric generators and that $\pi_1(L_b-\bfami_L)$ is generated
by a geometric basis $r_b$ then there is a presentation
$$
\pi_1(\CC^N-\hat\afami-\dfami)
\quad\cong\quad
\langle t_i,\hat r_a | t_i\inv\phi_b(t_i),\hat r_a t_i\inv\hat r_a\inv\phi_a(t_i),\rfami_a\rangle.
$$
where $\phi_a$ ($\phi_b$) is the automorphism associated to $r_a$ ($r_b$), $t_i$ is a free
geometric basis of $\pi_1(L-\dfami_L)$ and the $\hat r_a$ are lifts of $r_a$ by the topological section.
\end{lemma}

\proof
In the presentation of lemma \reff{pabd} we have simply
to set the geometric generators associated to $\hat\bfami$ to be trivial and to discard them from
the set of generators.
\qed

\begin{lemma}
\labell{pd}  
The fundamental group of $\CC^N-\dfami$ has a presentation
$$
\pi_1\quad\cong\quad
\langle t_i | t_i\inv\phi_b(t_i),\rho_at_i\inv \rho_a\inv\phi_a(t_i)\rangle.
$$
where
\begin{enumerate}
\item
the $t_i$ form a geometric basis of $\pi_1(L-\dfami_L)$,
\item
the $\rho_a$ can be expressed in terms of $t_i$ such that $\hat r_a\rho_a\inv$ is a geometric element
associated to $\hat\afami$ and a lift of $r_a$,
\item
the $\phi_a$, resp. $\phi_b$, are the braid monodromies associated to $r_a$, resp. $r_b$.
\end{enumerate}
\end{lemma}

\proof
We deduce this lemma using the presentation of lemma $3.6$.
Since each $\hat r_a$ is transversal to $\hat\afami$ it must be equal to a geometric element $\hat r_a'$ for
$\hat\afami$ up to some $\hat\rho_a$ expressible in terms of geometric elements for $\dfami$.
These elements in turn are expressible in terms of $t_i$ since the subgroup generated by
the $t_i$ is normal.

In the absence of $\hat\afami$ the geometric elements $\hat r_a'$ are obviously trivial, so from the presentation
of lemma $3.6$ we discard the generators $\hat r_a$ and we replace each $\hat r_a$ by $\hat\rho_a$
in the relations.
The relations $\rfami_a$ are thus replaced by a set of relations $\rfami(\hat\rho_a)$.
But these are in fact relations which we want to show to be superfluous.

So we choose $L''\subset H''\subset \CC^{N-1}$ such that
\begin{enumerate}
\item 
the preimages $E\subset H\subset\CC^N$ contain the line $L$,
\item
$L''$ and $H''$ are generic for $\afami$, thus in particular
$$
\pi_1(L''-\afami_{L''})\surj\pi_1(H''-\afami_{H''})\cong\pi_1(\CC^{N-1}-\afami).
$$
\end{enumerate}
These choices give rise to commutative diagrams
$$
\begin{array}{ccccccc}
\pi_1(L-\dfami_L) & \tto & \pi_1(E-\dfami_E-\hat{\afami}_E) &
\tto & \pi_1(L''-\afami_{L''}) & \tto & 1\\
\| & & \downarrow & & \downarrow\\[-4mm]
\| & & \downarrow & & \downarrow\\
\pi_1(L-\dfami_L) & \tto & \pi_1(H-\dfami_H-\hat{\afami}_H) &
\tto & \pi_1(H''-\afami_{H''}) & \tto & 1
\end{array}
$$
and
$$
\begin{array}{ccccc}
\pi_1(E-\dfami_E-\hat{\afami}_E) & \tto & \pi_1(E-\dfami_E) & \tto & 1\\
 \downarrow & & \downarrow\\[-4mm]
 \downarrow & & \downarrow\\
\pi_1(H-\dfami_H-\hat{\afami}_H) & \tto & \pi_1(H-\dfami_H) & \tto & 1.
\end{array}
$$
And we have the following presentations:
\begin{eqnarray*}
\pi_1(E-\dfami_E) & \cong &
\langle t_i|t_i\inv\phi_b(t_i),\rho_at_i\inv \rho_a\inv\phi_a(t_i)\rangle,\\
\pi_1(H-\dfami_{H}) & \cong &
\langle t_i|t_i\inv\phi_b(t_i),\rho_at_i\inv \rho_a\inv\phi_a(t_i), \rfami(\hat\rho_a)\rangle
\end{eqnarray*}
Our aim is to show that $E$ is sufficiently generic to imply 
$\pi_1(E-\dfami_E)\cong \pi_1(\CC^N-\dfami)$.

One way to see this is to proceed as follows:
Each relation in $\rfami(\hat\rho_a)$ is represented by a path in $E-\dfami_E$.
By our choice of $H$ for each such path there exists a $2$-cell in $H-\dfami_H$
with boundary freely homotopic in $H-\dfami_H$ to that path.

We will ultimately show that then each such path must be null-homotopic in $E-\dfami_E$
already, thus proving our claim.

But first we remark that each relation in $\rfami_a$ among elements $r_a$ in $\pi_1(H''-\afami_{H''})$
originates
in a singularity of the affine plane curve $\afami_{H''}$ or one of its asymptotes parallel to $L''$.

By genericity each singularity is either a cusp or a node and we may the $2$-cell we need in an arbitrarily close $3$-sphere
avoiding $\afami_{H''}$ and the line parallel to $L''$.

At the asymptotes parallel to $L''$ each relation is imposed by a singularity at infinity which is of $A$ type
and the $2$-cell can be found in a $3$-sphere around that point again avoiding $\afami_{H''}$ and the asymptote
parallel to $L''$ but also the line at infinity.

Each path and each $2$-cell is lifted to $H-\dfami_{H}$ via our topological section, so they correspond
precisely to the relations in $\hat\rfami(\hat\rho_a)$.
\\

Using projection of $H''$ along $L''$ to some affine line $\CC$, the total space $H-\dfami_{H}$
is mapped to $\CC$, such that outside a finite set in the target containing the values of singularities of $\afami_{H''}$
and of vertical asymptotes
we get a locally trivial fibration with fibre equal to $E-\dfami_E$.

By construction each $2$-cell belongs to the total space of this fibration.
So each path in the fibre $E-\dfami_E$ is null-homotopic in the total space due to the existence of the corresponding $2$-cell,
hence it must be null-homotopic in the fibre by the homotopy exact sequence 
taking into account that $\pi_2$ of a punctured affine line is trivial.
\qed

The claim of the lemma is of course only an intermediate step on our way to give a presentation
of the fundamental group.
Obviously we have to make the relations explicit in the sense that every relation is given in terms
of the chosen generators only.

Moreover we should try to reduce the number of relations whenever it is sensible to do so.

\paragraph{Remark}
We are very lax about the base points. They should be chosen in such a way
that all maps of topological spaces are in fact maps of pointed spaces.
(In particular in the presence of a topological section there is no choice left;
in the fibre and in the total space the base point is the intersection of the
section with the fibre and its projection to the base yields the base point there.)

\section{Brieskorn Pham unfolding}
\label{affine}

In this section we first construct a distinguished set of generators
for $\pi_1(\CC^N-\dfami)$.
We pick some distinguished fibres $L_v$ of the projection
$p_z:\CC^N\to\CC^{N-1}$ along the variable $z$ where in each case $L_v-\dfami_L$ can 
be equipped with a distinguished geometric basis by the method of Hefez and Lazzeri \cite{hl}.
For later use in section \reff{asympsection} we establish a relation between different such bases.

Secondly we want to give explicitly an exhaustive set of relations
associated to a geometric basis for the complement
of $\bfami$.
So we exploit the relation of two natural spaces
of perturbations of the Brieskorn-Pham polynomial
$f=y^3+x_1^\dd+\cdots+ x_n^\dd$.
On one hand there is our space $\CC^N$ considered as
the affine subspace of $\vnd$ of polynomials of weighted degree
at most $\dd$, which are monic as polynomials in $y$.
On the other hand there is a space germ $\CC^\mu,0$ of dimension
$\mu=2(d-1)^n$, the base of a versal unfolding of $f$ capturing the
perturbations of $f$ in the set-up of singularity theory.
That space and the fundamental group of its discriminant complement
has been under scrutiny in \cite{habil} and we will transfer some results to the current situation.

\subsection{Hefez Lazzeri path system}
\label{hl-path}

First we want to describe a natural geometric basis for some fibres of the projection
$p:\CC^N\to\CC^{N-1}$.
Since we follow Hefez and Lazzeri \cite{hl} we will call such bases accordingly.
We note first that fibres $L_u$ of the projection correspond to affine pencil of polynomials
$$
f_u(x_1,...,x_n)z
$$
and their discriminant points $\dfami_L$ are exactly the $z$ such that 
the $z$-level of $f$ is singular.
As in \cite{hl} we restrict our attention to the linearly perturbed
Brieskorn-Pham polynomial:
$$
f\quad=\quad y^3-3v_0y+\sum_{\ncount=1}^n (x_\ncount^\dd-\dd
v_\ncount x_\ncount).
$$

In that family the discriminant points for any generic pencil are in bijection
to the elements in the multiindex set of cardinality $2(\dd-1)^n$:
$$
\ind\quad=\quad\{\,(i_0,i_1,...,i_n)\,|\,
1\leq i_0<3,\,1\leq i_\nu< \dd-1\,\}.
$$
More precisely we get an expression for the critical values from 
\cite{hl}:

\begin{lemma}[Hefez Lazzeri]
\labell{laz-dis}
The polynomial defining the critical value divisor is given by the
expansion of the formal product ($\eta$ primitive of order $3d-1$)
$$
\prod_{\ibold \in\ind}
\left(-z+2(-1)^{i_0}v_0^\frac{3}{2}+
(\dd-1)\sum_{\ncount=1}^n
\eta^{i_\ncount}
v_\ncount^{\frac{\dd}{\dd-1}}\right).
$$
\end{lemma}

We deduce two immediate corollaries, that the discriminant sets
are equal for suitably related parameter values and that they
can be constructed inductively:

\begin{lemma}
\labell{twist}
The discriminant of the linearly perturbed polynomial $f$
is invariant under the multiplication of $v_0$ by a third
root of unity and of any $v_\ncount$ by a
$\dd$-th root of unity.
\end{lemma}

\proof
From the expansion above we see that the discriminant polynomial is a polynomial in
$v_\ncount^{\frac{\dd}{\dd-1}}$ but of course it is also a polynomial in 
$v_\ncount$, hence it must
be a polynomial in $v^\dd_\ncount$ since that is the least common 
power of both.
Then it is obviously invariant under multiplying $v_\ncount$ by a $d$-th root.
The statement for $v_0$ is proved analogously.
\qed

\begin{lemma}
\labell{circ}
The critical values of $f$ are distributed on circles of radius
$(\dd-1)|v_n|^{\frac{\dd}{\dd-1}}$
centred around the critical values of the polynomial
$$
f'\quad=\quad y^3-3v_0y+\sum_{\ncount=1}^{n-1}(x_\ncount^\dd
-\dd v_\ncount x_\ncount).
$$
\end{lemma}

\proof
Again we can use lemma \reff{laz-dis}.
A formal zero of the discriminant polynomial for $f$ differs by a term
$(\dd-1)v_n^{\frac{\dd}{\dd-1}}$ from a zero of the discriminant 
polynomial of $f'$, and that
difference is of the claimed modulus.
\qed

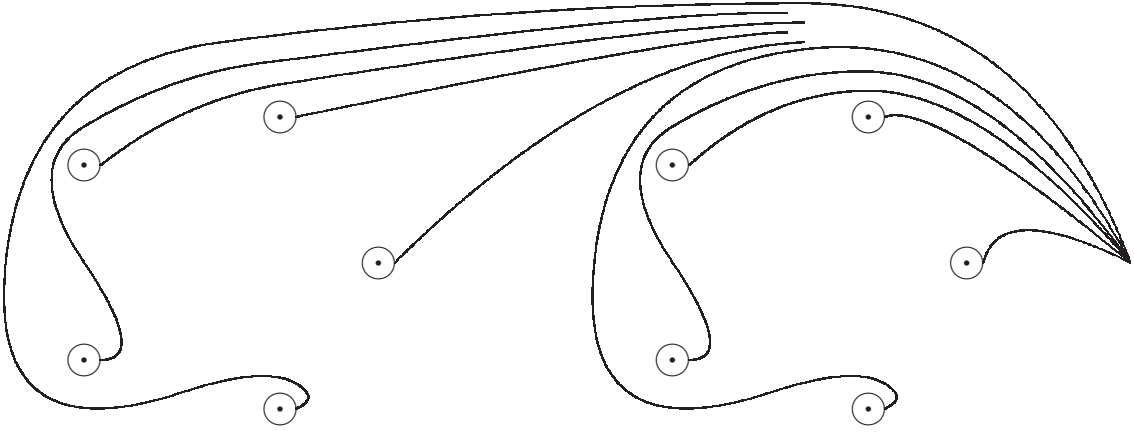
\begin{figure}[ht]
\begin{center}
\setlength{\unitlength}{.43mm}
\begin{picture}(180,150)(15,-70)

\put(50,0){\circle{10}}
\put(20,45){\circle{10}}
\put(-40,30){\circle{10}}
\put(-40,-30){\circle{10}}
\put(20,-45){\circle{10}}
\put(50,0){\circle*{2}}
\put(20,45){\circle*{2}}
\put(-40,30){\circle*{2}}
\put(-40,-30){\circle*{2}}
\put(20,-45){\circle*{2}}

\bezier{480}(55,0)(120,65)(180,68)

\bezier{480}(25,45)(150,71)(175,71)

\bezier{280}(-35,30)(-10,50)(20,55)
\bezier{480}(20,55)(140,74)(180,74)

\bezier{180}(-35,-30)(-20,-30)(-40,0)
\bezier{180}(-40,0)(-60,30)(-40,42)
\bezier{280}(-40,42)(-12,59)(18,62)
\bezier{480}(18,62)(140,78)(175,77)

\bezier{80}(25,-45)(30,-42.5)(28,-40)
\bezier{120}(28,-40)(20,-30)(-10,-40)
\bezier{280}(-10,-40)(-70,-60)(-64,0)
\bezier{280}(-64,0)(-58,60)(2,68)
\bezier{480}(2,68)(92,80)(180,80)
\bezier{280}(180,80)(250,80)(280,0)

\put(180,0){\begin{picture}(100,0)
\put(50,0){\circle{10}}
\put(20,45){\circle{10}}
\put(-40,30){\circle{10}}
\put(-40,-30){\circle{10}}
\put(20,-45){\circle{10}}
\put(50,0){\circle*{2}}
\put(20,45){\circle*{2}}
\put(-40,30){\circle*{2}}
\put(-40,-30){\circle*{2}}
\put(20,-45){\circle*{2}}

\bezier{130}(55,0)(60,20)(100,0)

\bezier{300}(25,45)(40,50)(100,0)

\bezier{280}(-35,30)(-10,53)(20,53)
\bezier{280}(20,53)(55,53)(100,0)

\bezier{180}(-35,-30)(-20,-30)(-40,0)
\bezier{180}(-40,0)(-60,30)(-40,42)
\bezier{280}(-40,42)(-12,59)(18,59)
\bezier{280}(18,59)(58,59)(100,0)

\bezier{80}(25,-45)(30,-42.5)(28,-40)
\bezier{120}(28,-40)(20,-30)(-10,-40)
\bezier{280}(-10,-40)(-70,-60)(-64,0)
\bezier{280}(-64,0)(-58,60)(2,66)
\bezier{280}(2,66)(62,72)(100,0)
\end{picture}}
\end{picture}
\caption{Hefez Lazzeri system in case $n=1$, $\dd-1=5$
\label{base}}
\end{center}
\end{figure}

We assume now that all $v_\ncount$ are positive real and of sufficiently distinct modulus
$$
|v_n|\ll\cdots\ll|v_1|\ll|v_0|.
$$
In case $n=1$ we define the Hefez Lazzeri geometric basis as indicated in figure \reff{base} for $\dd-1=5$,
where each geometric generator is depicted as a tail and a loop around a critical value.

Of course the geometric element associated to a loop-tail pair is represented by a closed
path based at the free end of the tail which proceeds along the tail, counterclockwise around
the loop and back along the tail again.

The $\dd-1$ elements of the base associated to punctures on the right
are denoted by $t_{1,1},...,t_{1,\dd-1}$, 
such that the corresponding critical values are enumerated 
counterclockwise starting on the positive real line.
Similarly the remaining elements are denoted by
$t_{2,1},...,t_{2,\dd-1}$.
Each tail on the right is thus labelled unambiguously by its second
index in $\{1,...,\dd-1\}$.

For the inductive step we suppose that the elements of a Hefez Lazzeri base for
$$
f'\quad=\quad y^3-3v_0y+
\sum_{\ncount=1}^{n-1}(x_\ncount^\dd-\dd
v_\ncount x_\ncount).
$$
are given by tail loop pairs each labeled by some multi-index in
$I_{n-1,d}$.
By assumption $v_n$ is sufficiently small compared to $v_{n-1}$ so 
we may assume that all critical values of
$$
f \quad=\quad y^3-3v_0y+
\sum_\ncount (x_\ncount^\dd-\dd v_\ncount x_\ncount).
$$
are inside the loops and in fact distributed at distance
$(\dd-1)|v_\ncount|^{\frac{\dd}{\dd-1}}$
from their centres.

In the inductive step each loop and its interior are erased
and replaced  by a scaled copy of the right hand side
of the Hefez Lazzeri base in the $n=1$ case, cf.\ figure \ref{base2}.
Each tail-loop pair with label $i_n\in\{1,...,\dd-1\}$ in an inserted disc 
fits with a tail labeled
by some $\ibold'=i_0i_1\cdots i_{n-1}$ to form a tail-loop pair 
representing an element of the new Hefez-Lazzeri base which is
labeled by $\ibold =i_0i_1\cdots i_{n-1}i_n$.

\newcommand{\fivedot}{\setlength{\unitlength}{.015mm}
\begin{picture}(100,0)(280,0)
\put(50,0){\circle{10}}
\put(20,45){\circle{10}}
\put(-40,30){\circle{10}}
\put(-40,-30){\circle{10}}
\put(20,-45){\circle{10}}
\put(50,0){\circle*{2}}
\put(20,45){\circle*{2}}
\put(-40,30){\circle*{2}}
\put(-40,-30){\circle*{2}}
\put(20,-45){\circle*{2}}
\end{picture}}

\begin{figure}[ht]
\begin{center}
\setlength{\unitlength}{.28mm}
\begin{picture}(280,150)(-110,-70)

\put(50,0){\circle{10}\fivedot}
\put(20,45){\circle{10}\fivedot}
\put(-40,30){\circle{10}\fivedot}
\put(-40,-30){\circle{10}\fivedot}
\put(20,-45){\circle{10}\fivedot}

\bezier{480}(55,0)(120,65)(180,68)

\bezier{480}(25,45)(150,71)(175,71)

\bezier{280}(-35,30)(-10,50)(20,55)
\bezier{480}(20,55)(140,74)(180,74)

\bezier{180}(-35,-30)(-20,-30)(-40,0)
\bezier{180}(-40,0)(-60,30)(-40,42)
\bezier{280}(-40,42)(-12,59)(18,62)
\bezier{480}(18,62)(140,78)(175,77)

\bezier{80}(25,-45)(30,-42.5)(28,-40)
\bezier{120}(28,-40)(20,-30)(-10,-40)
\bezier{280}(-10,-40)(-70,-60)(-64,0)
\bezier{280}(-64,0)(-58,60)(2,68)
\bezier{480}(2,68)(92,80)(180,80)
\bezier{280}(180,80)(250,80)(280,0)

\put(180,0){\begin{picture}(100,0)
\put(50,0){\circle{10}\fivedot}
\put(20,45){\circle{10}\fivedot}
\put(-40,30){\circle{10}\fivedot}
\put(-40,-30){\circle{10}\fivedot}
\put(20,-45){\circle{10}\fivedot}

\bezier{130}(55,0)(60,20)(100,0)

\bezier{300}(25,45)(40,50)(100,0)

\bezier{280}(-35,30)(-10,53)(20,53)
\bezier{280}(20,53)(55,53)(100,0)

\bezier{180}(-35,-30)(-20,-30)(-40,0)
\bezier{180}(-40,0)(-60,30)(-40,42)
\bezier{280}(-40,42)(-12,59)(18,59)
\bezier{280}(18,59)(58,59)(100,0)

\bezier{80}(25,-45)(30,-42.5)(28,-40)
\bezier{120}(28,-40)(20,-30)(-10,-40)
\bezier{280}(-10,-40)(-70,-60)(-64,0)
\bezier{280}(-64,0)(-58,60)(2,66)
\bezier{280}(2,66)(62,72)(100,0)
\end{picture}}

\put(-180,0){
\setlength{\unitlength}{.2mm}
\begin{picture}(100,0)
\put(0,0){\begin{picture}(200,200)(100,100)
\bezier{30}(0,100)(0,150)(40,180)
\bezier{30}(40,180)(100,220)(160,180)
\bezier{30}(160,180)(200,150)(200,100)
\bezier{30}(200,100)(200,50)(160,20)
\bezier{30}(160,20)(100,-20)(40,20)
\bezier{30}(40,20)(0,50)(0,100)

\bezier{60}(170,172)(227,118)(284,64)
\bezier{60}(130,4)(206,27)(282,50)
\end{picture}}

\put(50,0){\circle{10}}
\put(20,45){\circle{10}}
\put(-40,30){\circle{10}}
\put(-40,-30){\circle{10}}
\put(20,-45){\circle{10}}
\put(50,0){\circle*{2}}
\put(20,45){\circle*{2}}
\put(-40,30){\circle*{2}}
\put(-40,-30){\circle*{2}}
\put(20,-45){\circle*{2}}

\bezier{130}(55,0)(60,20)(100,0)

\bezier{300}(25,45)(40,50)(100,0)

\bezier{280}(-35,30)(-10,53)(20,53)
\bezier{280}(20,53)(55,53)(100,0)

\bezier{180}(-35,-30)(-20,-30)(-40,0)
\bezier{180}(-40,0)(-60,30)(-40,42)
\bezier{280}(-40,42)(-12,59)(18,59)
\bezier{280}(18,59)(58,59)(100,0)

\bezier{80}(25,-45)(30,-42.5)(28,-40)
\bezier{120}(28,-40)(20,-30)(-10,-40)
\bezier{280}(-10,-40)(-70,-60)(-64,0)
\bezier{280}(-64,0)(-58,60)(2,66)
\bezier{280}(2,66)(62,72)(100,0)
\end{picture}}

\end{picture}
\caption{Hefez Lazzeri system in case $n=2$, $\dd-1=5$
\label{base2}}
\end{center}
\end{figure}
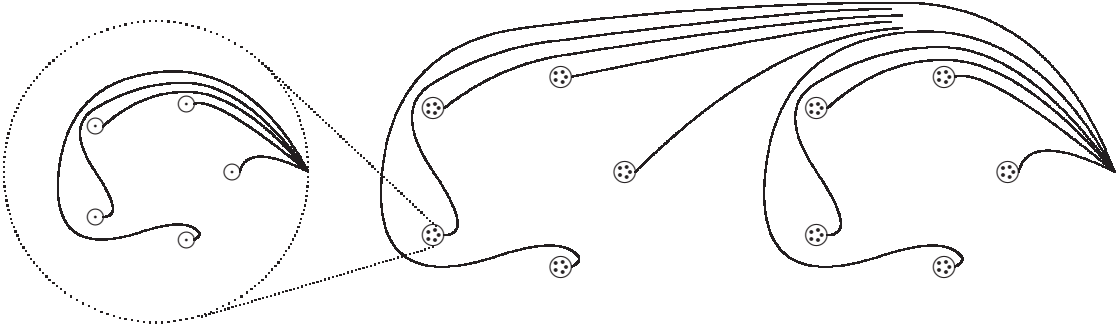

The element $\d_0$ represented by a path enclosing all critical values
counter-clockwise deserves our special attention.
In the case $n=1$ we see immediately that as an element in the 
fundamental group it can be expressed as
$$
t_{2,\dd-1}t_{2,\dd-2}\cdots t_{2,2}t_{2,1}
t_{1,\dd-1}t_{1,\dd-2}\cdots t_{1,2}t_{1,1}.
$$

To exploit the inductive construction for the general case we define
$t_{\ibold'}^+$,
$\ibold'\in I_{n-1,d}$, in the Hefez Lazzeri fibre of 
$$
f \quad=\quad y^3-3v_0y+
\sum_\ncount (x_\ncount^\dd-\dd v_\ncount x_\ncount)
$$
to be represented by the same tail-loop pair as the element $t_{\ibold'}$ in
the Hefez Lazzeri fibre of
$$
f'\quad=\quad y^3-3v_0y+
\sum_{\ncount=1}^{n-1}(x_\ncount^\dd-\dd v_\ncount x_\ncount).
$$

\begin{lemma}
\labell{bundle}
For any $\ibold'\in I_{n-1,d}$ there is a relation
$$
t_{\ibold'}^+\quad=\quad t_{\ibold'(\dd-1)}t_{\ibold'(\dd-2)}
\cdots t_{\ibold'2}t_{\ibold'1},
$$
where ${\ibold'(\dd-1)},{\ibold'(\dd-2)},\ldots, {\ibold'2},{\ibold'1}$
are the obvious elements of $\ind$.
\end{lemma}

\proof
The loop of $t_{\ibold'}^+$ is the path which encloses
counter-clockwise $\dd-1$ of the
critical values which are inserted in the inductive step.
So it is homotopic to the product with descending indices of the
loop-tail pairs inserted into that loop.

This relation is preserved under appending the tail of $t_{\ibold'}$
and so the claim follows.
\qed

To formulate the general claim we define a linear order $\prec_0$ on
$\ind$ to be
the lexicographical order with respect to the linear order $>$ (!)
on each component.
Accordingly we introduce an order preserving enumeration function
$$
\ibold_0:\,(\{1,...,2(\dd-1)^n\},<)\,\tto\,(\ind,\prec_0).
$$

\begin{lemma}
\labell{cox}
Suppose the element $\d_0$ is represented by a path which
encloses counter-clockwise
all critical values of a Hefez Lazzeri fibre, then
$$
\d_0\quad=\quad\prod_{k=1}^{2(\dd-1)^n}t_{\ibold_0(k)}.
$$
\end{lemma}

\proof
In case $n=1$ we have given above an expression for $\d_0$ which is the expression claimed
here, as we have taken care that the order given by $\prec_0$ on $I_{1,d}$ is
$$
2,\dd-1 \,\prec_0\, 2,\dd-2 \,\prec_0\, \cdots\,\prec_0\, 2,1
 \,\prec_0\,1,\dd-1 \,\prec_0\, 1,\dd-2 \,\prec_0\, \cdots\,\prec_0\, 1,1.
$$
Inductively we then get an expression for $\d_0$ from that
of $\d_0'$ using the definition of $t_{\ibold'}^+$:
$$
\d_0'\quad=\quad
\prod_{k=1}^{2(\dd-1)^{n-1}}t_{\ibold_o'(k)}.
\qquad\implies\qquad 
\d_0\quad=\quad
\prod_{k=1}^{2(\dd-1)^{n-1}}t_{\ibold_o'(k)}^+
$$
Finally it suffices to replace each $t_{\ibold'}^+$ using
lemma \reff{bundle} above to get our claim.
\qed

\labell{pathsandhomotopies}

As we noticed in lemma \reff{twist} the set of singular values
remains unchanged upon multiplication of any of the real
$v_\ncount$ by a $\dd$-th root of unity $\xi$,
resp.\ $v_0$ by a third root.
The corresponding fibres are thus equipped with the same 
Hefez Lazzeri systems of paths.

We denote by $t_\ibold(\jbold)$, $\jbold\in\ind$,
the elements of the Hefez-Lazzeri basis
in the fibre at
$v(\jbold)=(v_0\xi^{d(j_0-1)}v_1\xi^{j_1-1},...,v_n\xi^{j_n-1})$
and by $z_0$ the fibre coordinate of the Hefez Lazzeri base
point, which may be assumed to belong to a topological section
as in the
proof of lemma \ref{sabd}.

We can now compare the fundamental groups
$\pi_1(L_{v(\jbold)}-\dfami_{L_\jbold},(v(\jbold),z_0))$
along paths
$$
\oo_\jbold:\quad s\,\mapsto\,(v_0\xi^{sd(j_0-1)}v_1\xi^{s(j_1-1)},...,
v_n\xi^{s(j_n-1)},z_0).
$$

\begin{lemma}
\labell{hotopy}
Conjugation by a path $\oo_\jbold$ induces an isomorphism
$$
\oo^*_\jbold:\pi_1(L_{v(\jbold)}-\dfami_{L_\jbold},(v(\jbold),z_0))\to\pi_1(L_v-\dfami_L,(v,z_0))
$$
such that for $\ibold_0(1)=11\cdots 1\in\ind$ 
$$
\oo^*_\jbold(t_\jbold(\jbold))\quad=\quad t_{\ibold_0(1)}
$$
\end{lemma}

\proof
We consider the case $j_0=...=j_{n-1}=0$ first.
Then along $\oo_\jbold$ all punctures move counterclockwise
in the innermost inserted discs covering an angle of $(j_n-1)\varth$,
$\varth=\frac{2\pi}{\dd-1}$.

Accordingly the final part of each tail has to be adjusted alongside,
but all other parts may just be kept fixed.
In particular the tail segment with label $j_n$ is moved to the 
segment with label $1$.

Similarly for any $j_\ncount\neq1$, the $\ncount$-th segment
of each tail is affected.
While precursory segments are unaffected in this more general case
successive segments are moved as well, 
but they are moved to segments having the same label.

We may conclude that in the general case the segments labeled
by the components of $\jbold$ are moved to segments labeled by
$1$ and so our claim holds.
\qed

\subsection{Brieskorn-Pham monodromy}

The next aim is to determine the set of relations imposed on Hefez-Lazzeri generators
by the geometric generators associated to $\bfami$.

Our strategy is to use the relations imposed by the geometric elements associated to the
bifurcation set $\bfami_f$ in a truncated versal unfolding of the Brieskorn-Pham polynomial
$$
f=y^3+x_1^\dd+x_2^\dd+\cdots+ x_n^\dd,
$$
with base the affine space germ $\CC^{\mu-1},0$, where $\mu$ is the
Milnor number of $f$.

We can do so by means of the \emph{truncated subdiagonal unfolding}, the unfolding
of $f$ by non-constant monomials of degree less than $d$.
Its affine base $A$ with bifurcation set $\bfami_A$ is
naturally a subspace of $(\CC^{N-1},\bfami)$, the germ at $0$ a
subgerm of $(\CC^{\mu-1},\bfami_f,0)$.

In fact our aim is to produce a projection $p_v:\CC^{N-1}\to\CC^{N-2}$ as in lemma $3.4$
and a suitable fibre $L_v$ of it which belongs to $A$.

For the following discussion we observe that a fibre $L_v$ defines a 
pencil of polynomials, which are linearly related in the sense, 
that their differences are scalar multiples of a linear polynomial.
Moreover a critical point with value $z_0$ of any such polynomial
corresponds to a critical point on the hypersurface 
defined as its $z_0$-level.

We will require two conditions on $L_v$ which are specified in the following results:

\begin{lemma}
\labell{a2gen}
There is a pencil of linearly related polynomials in $A$, each of which has only non-degenerate critical points,
except for a finite number of polynomials with degenerate critical points of
type $A_2$ only.
\end{lemma}

\proof
First we consider the case $n=1$. In this case we can perturb
$x_1^\dd$  to a polynomial $f_1$ of degree
$\dd$ in $x_1$ such that the pencil $y^3+v_0y+f_1+v_1x_1$
of linear perturbations has the required property.

For the general case we just take the sum of such perturbations and thus get a family depending on $n$ parameters $v_k$, which
are the coefficients of the linear monomials:
$$
y^3+y+f_1+v_1x_1 +\cdots+ f_n+v_nx_n.
$$ 

Since the critical points of the sum are the points such that each coordinate is a critical point of the corresponding
summand and the Hessian is diagonal, 
we deduce, that to get a degenerate critical point at least in one coordinate the corresponding critical
point must be degenerate.

In fact it is non-generically degenerate if and only if either the critical point is degenerate in two coordinates
or it is non-generically degenerate in one coordinate.

But with our choice this happens only in a codimension two set.
Hence we can find a pencil as claimed in the $n$-dimensional family obtained from the one-dimensional
perturbation in each coordinate.
\qed

We remark that if a polynomial in $x_1$ has one degenerate critical point, then
its sum with a polynomial in $x_2$ has several degenerate critical
points in bijection to the
critical points of the second summand.

\begin{lemma}
\labell{a11gen}
There is a pencil of linearly related polynomials in $A$, such that
all polynomials have critical points with distinct values, except
for a finite set of polynomials, for which there is one multiple value.
It must
belong to a single pair of non-degenerate critical points.
\end{lemma}

\proof
We consider the case $n=1$.
If in the variable $x_1$ we pick a pencil $f_1+v_1x_1$ of polynomials 
as in the proof above,
there are at most two coinciding critical values for each polynomial.
Hence there is a lower bound $\varepsilon>0$ such that at most one pair of critical values has distance less than $\varepsilon$.

In the second variable we pick a generic pencil $f_2+v_2x_2$ of linearly related polynomials
with a polynomial $f_2$ of which all critical values and their
difference are distinct and within the bound $\varepsilon$
of each other.

Their sum yields a two-dimensional family of linearly related polynomials
$$
f_1+v_1x_1+f_2+v_2x_2.
$$
Its critical values are sums of the critical values of both summands.

We deduce that the subfamily $f_1+v_1x_1 +f_2$ has only polynomials with at most two coinciding critical values.

Suppose now one of the critical points involved were degenerate.
Then we can choose a sufficiently close polynomial such
that three critical values are arbitrarily close in contradiction to our
construction.

So the corresponding line -- and sufficiently close parallel lines --
in the parameter plane $v_1,v_2$ do not meet the locus of polynomials with non-generically
conflicting values. Thus this locus must be contained in codimension two, i.e.\ points, or in parallel lines.

Suppose there is such a parallel line corresponding to a family
$f_1+v_1x_1+f_2'$. Then there is a constant $a$, such that the constant
summand and all polynomials of the variable summand have
critical values of difference $a$. Now the set of critical values
of the family $f_1+v_1x_1$
is described a monic discriminant polynomial $d=d(v_1,z)$ which is
irreducible for a generic choice.
The condition on the distance of values translates into the condition
$$
\ON{res}_z(d(v_1,z), d(v_1,z-a)) = 0.
$$
But then the two polynomials have a common factor and are thus equal
by their irreducibility,
which is possible only for $a=0$. So such a line does not exists
and we have shown that polynomials with bad critical value
properties occur only in codimension $2$.
\\

In the general case $n\geq2$ we argue as before, but start with
a generic pencil in all but the last variable.
\qed

The full genericity property we need comprises the two properties above and the property, that projection along
$L'$ is generic in the sense of lemma $3.4$.

\begin{lemma}
\labell{b-gen}
There is a projection as in lemma $3.4$ such that a fibre $L_v$
corresponds to pencils as in lemma $4.7$ and $4.8$.
\end{lemma}

\proof
By lemma $3.4$ we know that there is a projection $p_v:\CC^{N-1}\to\CC^{N-2}$
such that a generic fibre is generic for $\bfami$. Notice that we have an open condition, hence a Zariski open
set of projections has that property.

Similarly the properties of the two preceding lemmas $4.7$ and
$4.8$ put open conditions on the choice of the projection.
Hence we conclude, that a projection $p_v$ and a fibre $L_v$ can be
found as claimed.
\qed

We remark that it is not true, that $L_v$ necessarily is a generic fibre
of $p_v$ but still it is sufficiently good.
We use $L_v$ and an arbitrarily close generic parallel
$L'$ to compare the braid monodromy relations of \cite{habil}
with the relations imposed by geometric elements associated to $\bfami$ in
$\pi_1(\CC^{N}-\dfami)$.
To facilitate the use of citations from \cite{habil}
let us recall:

\begin{quote}
Given a linear projection $\CC^2\to\CC$ which induces a finite map
of a plane curve $C\subset\CC^2$ to $\CC$, any closed path in the
image avoiding the bifurcation set of $C$ gives rise to a braid and the
associated isotopy class of diffeomorphisms of the punctured fibre.

The braids thus obtained form the \emph{braid monodromy
group}. It acts naturally on the free group generated by a geometric 
basis in the punctured fibre.

Generalisations to divisors in affine spaces or to the analogue situation
for germs are immediate.
\end{quote}

To pursue our argument further let $t_i$ be a Hefez Lazzeri geometric base for
the fibre $L$ of the projection $p:\CC^N\to\CC^{N-1}$ over a point corresponding
to a linear perturbation of the polynomial $f$.

We may assume that this polynomial belongs to the pencil determined by $L_v$.

\begin{prop}
\labell{kernel}
Suppose $\pi_1(L'-\bfami)$ is generated by a geometric basis $r_b$ with base
point close to point to which $L$ projects. Suppose further that the braid monodromy
group $\Gamma$ of the discriminant in the base of
the semi-universal unfolding of $f$ is generated
by elements $\{\b_s\}$.
Then there is an identity of normally generated normal subgroups 
$$
\langle\langle t_i\inv\phi_b(t_i)\rangle\rangle\quad=\quad\langle\langle t_i\inv\b_s(t_i)\rangle\rangle
$$
\end{prop}

\proof
Arbitrarily close to $L_v$ there is a parallel line $L'_f$ in the base $\CC^{\mu-1}$
of the truncated unfolding of $f$, which is generic with respect to $\bfami_f$.
We denote now by $E',E_v$ and $E_f$
the complex affine planes which project to
the lines $L',L_v$ and $L'_f$ respectively. By intersection with the
respective discriminants they contain discriminant curves $\dfami_{E'},\dfami_{E_v}$
and $\dfami_{E_f}$.

By the result of \cite{habil} the normal subgroup on the right hand side
is the kernel of the natural map from the free group generated by the $t_i$
to $\pi_1(E_f-\dfami_{E_f})\cong\pi_1(\CC^\mu-\dfami_f)$. By the genericity
of $L'$ the left hand side is the kernel of the map to $\pi_1(E'-\dfami_{E'})$.

Hence both sides are equal, if we show, that the curves $\dfami_{E_f}$ and
$\dfami_{E'}$ are isotopic through a path of plane curves.
\\

For the isotopy claim we first observe that both curve naturally degenerate to 
$\dfami_{E_v}$.
By lemma $4.7$ and $4.8$ above $\dfami_{E_v}$ has only double points.
Since parallels to $L_v$ are transversal to $\bfami_f$ in $\CC^{\mu-1}$
the degeneration from $\dfami_{E_f}$ to $\dfami_{E_v}$ is equisingular and so they are
isotopic (though not via an isotopy which preserves the projection).
Hence $\dfami_{E_v}$ is -- like $\dfami_{E_f}$ -- a plane curve with only ordinary cusps
and ordinary nodes.
Therefore the degeneration from $\dfami_{E'}$ to $\dfami_{E_v}$ is also equisingular,
so $\dfami_{E'}$ and $\dfami_{E_f}$ are isotopic.
\qed

The explicit description of the normal subgroup in term of a finite set of relations will be taken from \cite{habil}
and used in the proof of the main theorem in section \ref{results}.

\section{Asymptotes}
\label{asympsection}

Having fixed a preferred choice of generators for the discriminant complement in the previous section,
it is now time to get the explicit relations imposed by the degenerations along the divisor $\afami$.

The first and easier task is to understand a local model, that is
the restriction of the projection along $z$ to the preimage of
a small disc transversal
to $\afami$.
In fact we will get a presentation for its fundamental
group in terms of a geometric basis in a local reference fibre
subjected to a single relation.

Next we address the hard problem -- which lies at the heart
of our argument --  to
bring those local relations together globally:
Basically we can start with a path connecting a global reference fibre
to a local disc transversal to $\afami$. The trivialisation along the path
induces an isomorphism between the fundamental groups of 
the local and global reference fibres at the endpoints.

In practise we are just able to determine unambiguously the images
of local
relations along a very restricted set of paths, which we
have to construct with great care.

Finally we notice that a geometric element associated to $\afami$
comes naturally with each path and the boundary of the corresponding
local disc.
In fact it is our final task for this section to show that the geometric
elements thus obtained suffice to apply the results of section $3$.

\subsection{affine and local models}

\newcommand{\order}{n}
\newcommand{\Cpunc}{\CC-{\scriptstyle[}\order{\scriptstyle]}}  
\newcommand{\Cppunc}{\CC-[2]}  
\newcommand{\Dpunc}{D^0-{\scriptstyle[}m{\scriptstyle]}}  

To describe the geometry of the projected discriminant locally at a point of
$\afami$ we consider first the plane affine curves $C_\order$
\begin{eqnarray*}
C_\order\subset \CC^2 & : & y(yx^\order-1) \quad=\quad 0.
\end{eqnarray*}

We call $C_\order$ an \emph{asymptotic curve} of order $\order$,
since its equation can be rewritten as $y(y-\frac1{x^\order})$ to show that
the $y$-axis is a vertical asymptote and $\order$ is its pole order.

\begin{lemma}
\labell{asyp-I}
Suppose $b_0,b$ form a geometric basis in a vertical fibre
associated to $y=0$ and $yx^n=1$ respectively. Then the complement
of the curve $C_n$ in $\CC^2$ has fundamental group presentable as
\begin{eqnarray*}
\pi_1(\CC^2-C_{n}) & \cong & \langle b_0,b\,|\, b^nb_0 = b_0b^n \rangle
\end{eqnarray*}
\end{lemma}

\proof
The complement of $C_{n}$ is isomorphic to the complement of $C'_{n}$, 
a plane affine curve given by
$$
z(x^n-z),
$$
since both are complements of a quasi-projective curve in
$\CC\times\PP^1$ given by $yz(x^ny-z)$.

The projection along $z$ exhibits the complement of $C'_{n}$ to be a fibre
bundle over the punctured complex line $\CC^*$.
The fibre is diffeomorphic to a $2$-punctured complex line, 
which we denote by $\Cppunc$.
Since the projection map has a section, we get a short exact sequence 
with split surjection
$$
1\to \pi_1(\Cppunc)\to\pi_1(\CC^*\times\CC-C'_{n})
\to \pi_1(\CC^*)\to 1.
$$
With base points $z_0>1,x_0=1$ we may choose elements
\begin{itemize}
\item
$a_0,b_0$ in $\pi_1(\Cppunc)$ represented by a geometric elements
associated to $z=0$ and $z=x^n$ respectively,
\item
$c$ in $\pi_1(\CC^*\times\CC-C'_{n})$ represented by a geometric element in the
line $z=z_0$ associated to $x=0$.
\end{itemize}

All relations in $\pi_1(\CC^2-C'_{n})$ follow from the conjugacy action
of $c$ on $a_0,b_0$ in
$\pi_1(\CC\times\CC^*-C'_{n})$, which can be determined as a braid
monodromy,
and the triviality of $c$ in $\pi_1(\CC^2-C'_{n})$.

Indeed our claim is derived from the conjugacy action given explicitly as
$$
c a_0 c\inv = (a_0b_0)^n a_0(a_0b_0)^{-n},\quad
c b_0 c\inv = (a_0b_0)^n b_0(a_0b_0)^{-n},
$$
which is an immediate consequence of the braid monodromy given by $n$ full twists $\s_1^{2n}$.
Since either $a_0b_0b=1$ or $b_0a_0b=1$, we should replace 
$a_0b_0$ by $b\inv$ or $b_0bb_0\inv$ and we get the claim in either case.
\qed

\paragraph{Remark:}
Note that the complement of $C_\order$ can be strongly retracted to any
fibred neighbourhood of the $y$-axis,
since $C_\order$ is invariant under $(x,y)\mapsto (\l^\order y,\l\inv x)$.

\paragraph{Definition}
\begin{sloppypar}
A curve $C\subset D\times\CC$ is called an \emph{asymptotic germ}
of fibre degree $m$ and pole order $n$ if
\begin{enumerate}
\item
the fibration to $D$ is locally trivial onto the punctured disc $D^*$
with fibre a \mbox{$m$-punctured} complex line, which we denote by
$\CC-{\scriptstyle[}m{\scriptstyle]}$,
\item
the projective closure of $C$ in $D\times\PP^1$ is a topological cover
of $D$ of degree $m$ and intersects the line at infinity with multiplicity
$\order$.
\end{enumerate}
\end{sloppypar}

An example with fibre degree $m+1$ and pole order $\order$ is
provided by the plane curve 
$$
(y^m-1)(x^ny-1).
$$
\\

A neighbourhood of the $y$-axis in $\CC\times\PP^1$ is 
covered by two bi-disc
$D\times D^0$ and $D\times D^\infty$ centered at $0$ and $\infty$,
which intersect in the product of $D$ with an annulus $A$.
Given any asymptotic germ $C_{n,m}$ of fibre degree $m+1$ and pole
order $\order$, that kind of open cover can be chosen to cover
a sufficiently small neighbourhood $D\times\PP^1$ such that
\begin{itemize}
\item
 $D\times A$ and $C_{n,m}$ are disjoint,
\item
$D\times D^\infty$ contains only the branch of $C_{n,m}$ through $\infty$, 
\item
$D\times D^0$ contains all other branches of $C_{n,m}$.
\end{itemize}
For $D$ small enough $D\times D^0$ meets $C_{n,m}$ in $m$ disjoint
horizontal copies of $D$, thus we can identify the complement with the
product of $D$ with the
$m$-punctured disc $D^0$, which we denote by $\Dpunc$.

On the other hand $C_{n,m}$ meets $D\times D^\infty$
in a single smooth branch, which
intersects the line at infinity with multiplicity $n$. 
And so does the branch of our model curve $C_n$ in the complement
of the $x$-axis $y=0$.
Therefore we can discard the branch of $C_{n,m}$ and the line at
infinity from $D\times D^\infty$ to get a complement,
which can be identified with the complement of $C'_{n}$ in $D\times
\CC$.
Hence we get a decomposition
$$
(D\times\CC)-C_{n,m} = (D\times\CC-C_{n})\cup_{D\times A}
(D \times(\Dpunc)).
$$

\begin{prop}
\labell{local1}
If $C$ is an asymptotic germ of fibre degree $m+1$ and pole order $n$ then
there is a geometric basis $b_0,b_1,...,b_m$ in a vertical fibre of the complement such that
$$
\pi_1(D\times\CC-C) \quad \cong \quad \langle b_0,b_1,...,b_m\,|\,
(b_1\cdots b_m)^n b_0=b_0 (b_1\cdots b_m)^n\rangle
$$
with $b_0$ associated to the branch going to infinity and $D$
sufficiently small.
\end{prop}

\proof
By the hypothesis on $D$ we
can employ the decomposition discussed above.
We can choose a geometric basis such that
the $b_i$ belong to $D\times (\Dpunc)$ except for $b_0$, which
belongs to the part around infinity.
Then $\pi_1(D\times(\Dpunc)) $ is the free group generated freely by $b_1,...,b_m$
and $\pi_1(D\times \CC- C_n)$ is isomorphic to $\pi_1(\CC^2-C_n)$ and can by lemma \ref{asyp-I} be presented as
$$
\langle b_0,b \,|\,b^nb_0 = b_0 b^n \rangle.
$$
Here $b$ is the geometric generator associated to $y=0$, which
is to be identified with the core of the annulus $A$.

We can then apply the van Kampen theorem
to our decomposition with the observation that $b_1\cdots b_m$
and $b$ represent the same generator of the fundamental group of the
intersection,
and we get
\begin{eqnarray*}
\pi_1 & \cong &
\langle b_0,b_1,...,b_m,b \,|\, b^nb_0 = b_0 b^n,\, b=b_1\cdots b_m \rangle
\end{eqnarray*}
from which the claim is immediate by removing $b$.
\qed

A geometric basis as in the proposition is called \emph{adapted} to the
asymptote and provides a tool to make the second set of relations in lemma
\ref{pd} more explicit using the element $\d_0$ introduced in lemma
\ref{cox}.




\newcommand{\sumld}{\Sigma_{\l,d}}
\newcommand{\sumk}{\mbox{$\sum_\ncount$}}
\newcommand{\lk}{\l_\ncount}
\newcommand{\rk}{\rho_\ncount}
\newcommand{\xk}{x_\ncount}
\newcommand{\sumi}{\mbox{$\sum_\iota$}}
\newcommand{\li}{\l_\iota}
\newcommand{\ri}{\rho_\iota}
\newcommand{\rplus}{\RR^+}

\newcommand{\pt}{\l}
\newcommand{\rsol}{\hat x}
\newcommand{\csol}{\check x}

\newcommand{\lcrit}{\l_{\text{crit}}}
\newcommand{\ptcrit}{\pt_{\text{crit}}}

\newcommand{\vy}{\l_0}
\newcommand{\vi}{\l_i}
\newcommand{\vn}{\l_n}

\subsection{a sufficient subfamily}

In this subsection we consider a suitable family $\gfami$ of
polynomials in which our parallel transport will take place.
$$
\begin{array}{cccl}
f  & = & &
y^3-3\vy y + \sum_i \left( x_i^\dd-\dd \vi x_i \right)\\[2mm]
&&
+ & x_n^\dd-\dd \vn x_n-\frac{3d}{3d-1}\vn' x_n^{3d-1}
-  3\pt \left( \vy y x_n^{2d} + d\sum_i \vi x_i x_n^{\dd-1}\right)
\end{array}
$$
We use
the following claims to trace a distinguished
critical value along suitable paths.

\begin{lemma}
\labell{crit-val}
The value of $f$ at critical points is given by
$$
f|_{\nabla f=0} \quad=\quad
-2\vy y - (3d-1) \vn x_n -
{\textstyle\frac{1}{3d-1}}\vn'x_n^{3d-1}
 - (3d-1) \sum_i \vi x_i.
$$
\end{lemma}

\proof
First the vanishing gradient condition may be expressed by the
following system of equations:
\begin{eqnarray}
y^2 & = & \vy(1+\pt x_n^{2d})\notag
\label{crit-eqn}
\\
x_i^{3d-1} & = & \vi(1+\pt x_n^{3d-1})\\
x_n^{3d-1} & = & \vn + \vn' x_n^{3d-2}
+ \pt(2\vy y x_n^{2d-1}
+(3d-1) \sum_i \vi x_i x_n^{3d-2})\notag
\end{eqnarray}
Next we use them to replace the pure monomials in $f$:
\begin{eqnarray*}
f & = &
y v(1+\pt x_n^{2d}) -3v y 
+\sum_i \left( x_i(\vi-\pt d \vi x_n^{3d-1})-3 d \vi x_i \right) \\
&&
+x_n(\vn + \vn' x_n^{3d-2}+ 2\pt v y x_n^{2d-1} 
+(3d-1)\pt \vi x_i x_n^{3d-2})\\[4mm]
&&
-\dd \vn x_n-{\textstyle\frac{3d}{3d-1}}\vn' x_n^{3d-1}
-3\pt \big(v y x_n^{2d} +d \sum_i \vi x_i x_n^{3d-1} \big)\\
& = &
-2\vy y -(3d-1) \vn x_n
-{\textstyle\frac{1}{3d-1}}\vn'x_n^{3d-1}
 -(3d-1)  \sum_i \vi x_i.\\[-14mm]
\end{eqnarray*}
\qed

\begin{lemma}
\labell{simple-asym}
Suppose $\vy^3\in\rplus$, $\vi^{3d}\in\rplus$ and $\vn'\in\rplus$.
Then there is a positive real $\ptcrit$ such that the number of 
critical points (counted with multiplicity) is maximal for
$\pt\in[0,\ptcrit[$ and drops by one at $\pt=\ptcrit$.
\end{lemma}

\proof
From the equations (\ref{crit-eqn}) we deduce a polynomial equation
for the last coordinate $x_n$ of all critical points.
To eliminate the other coordinates we note, that the following
expression with $\xi_i,\xi_0$ primitiv roots of unity
of order $3d-1$ resp.\ $2$,
$$
\prod_{\rho=1}^2\prod_{\rho_i=1}^{3d-1}
\big(x_n^{3d-1}-\vn-\vn'x_n^{3d-2}- 2\pt \xi_0^\rho\, \vy y x_n^{2d-1}
- (3d-1)\pt \sum_i \xi_i^{\rho_i} \vi x_i x_n^{3d-2}
\big)
$$
is zero on critical points by the last equation of (\ref{crit-eqn}).
Due to the invariance under $x_i\mapsto \xi_i x_i$,
$y\mapsto \xi_0 y$, our expression is a polynomial $y^2$ and
$x_i^{3d-1}$.
We may thus insert the expressions given by the right hand sides of
(\ref{crit-eqn}) to get a polynomial in $x_n$ only.

To understand the leading and subleading coefficient, we may neglect
lower order terms in the elimination process. Then we get the formal
expression
$$
\prod_{\rho=1}^2\prod_{\rho_i=1}^{3d-1}
\big(x_n^{3d-1}-\vn'x_n^{3d-2}
- 2\pt^{\frac32}\xi_0^\rho\,  \vy^{\frac32} x_n^{3d-1}
-(3d-1)\pt^{\frac{3d}{3d-1}}\sum_i \xi_i^{\rho_i} 
 \vi^{\frac{3d}{3d-1}}x_n^{3d-1}
\big),
$$
where the $\vi$ are uniquely determined as positive real numbers
thanks to the hypotheses $\vy^3$, $\vi^{3d}\in \rplus$.
So the leading coefficient is given by
$$
\prod_{\rho=1}^2\prod_{\rho_i=1}^{3d-1}
\big(1
- 2\pt^{\frac32}\xi_0^\rho\, \vy^{\frac32}
- (3d-1)\pt^{\frac{3d}{3d-1}} \xi_i^{\rho_i}\, \vi^{\frac{3d}{3d-1}}
\big)
$$
Let $\ptcrit$ be its smallest positive real root in $\pt$, which
is in fact that of the factor with $\xi_0,\xi_i=1$.
Then the algebraic number of critical points drops at $\ptcrit$
for the first time.
Moreover the next coefficient is given by the sum of all the different
possibilities to neglect one factor of the product up to the constant
$\vn'>0$, hence it is non-zero at $\pt=\ptcrit$ if and only if no other
factor except the obvious one vanishes, which is obviously the case. 
\qed

\begin{lemma}
\labell{hessian}
Suppose $\vn', \vn\in\RR^{\geq0}$, $\vy,\vi\neq0$.
Then there is no degenerate critical point with
$\vy y$, $\vi x_i$, $x_n \in \rplus$, $\pt\in[0,\ptcrit[$,
except if all $\pt,\vn$ and $\vn'$ vanish.
\end{lemma}

\proof
It suffices to show that the gradient vectors to the equations
(\ref{crit-eqn})
are linearly independent at solutions with $\vy y$, $\vi x_i \in \rplus$.
\begin{eqnarray*}
\begin{vmatrix}
2y & 0 & \cdots & -2d \pt \vy x^{2d-1}\\
0 & (3d-1)x_i^{3d-2} & \cdots & -\pt(3d-1)\vi x_n^{3d-2}\\[1.5mm]
\vdots & \vdots & \ddots & \vdots\\[1.5mm]
-2\pt \vy x_n^{2d-1} & -(3d-1)\pt \vi x_n^{3d-2} & \cdots &
\begin{array}{c}
(3d-1)x_n^{3d-2} -(3d-2)\vn' x_n^{3d-3}\\
-2(2d-1) \pt \vy yx_n^{2d-2}\\
-(3d-1)(3d-2)\pt \sum_i \vi x_i x_n^{3d-3}
\end{array}
\end{vmatrix}
& \neq & 0\\
\end{eqnarray*}
We multiply each column by the corresponding variable, which is
non-vanishing by assumption to get an equivalent claim
\begin{eqnarray*}
\begin{vmatrix}
2y^2 & 0 & \cdots & -2d \pt \vy x^{2d}\\
0 & (3d-1)x_i^{3d-1} & \cdots & -\pt(3d-1)\vi x_n^{3d-1}\\[1.5mm]
\vdots & \vdots & \ddots & \vdots\\[1.5mm]
-2\pt \vy yx_n^{2d-1} & -(3d-1)\pt \vi x_i x_n^{3d-2} & \cdots &
\begin{array}{c}
(3d-1)x_n^{3d-1} -(3d-2)\vn' x_n^{3d-2}\\
-2(2d-1) \pt \vy yx_n^{2d-1}\\
-(3d-1)(3d-2)\pt \sum_i \vi x_i x_n^{3d-2}
\end{array}
\end{vmatrix}
& \neq & 0\\
\end{eqnarray*}
Next we apply the equations (\ref{crit-eqn})
and simplify the entry in the right bottom corner
\begin{eqnarray*}
\begin{vmatrix}
2\vy(1+\pt x_n^{2d}) & 0 & \cdots & -2d \pt \vy x^{2d}\\
0 & (3d-1)\vi(1+\pt x_n^{3d-1}) & \cdots & -\pt(3d-1)\vi x_n^{3d-1}
\\[2mm]
\vdots & \vdots & \ddots & \vdots\\
-2\pt \vy y x_n^{2d-1} & -(3d-1)\pt \vi x_i x_n^{3d-2} & \cdots &
\begin{array}{c}
(3d-1)\vn+\vn' x_n^{3d-2}\\
+2d \pt \vy yx_n^{2d-1}\\
+(3d-1)\pt \vi x_i x_n^{3d-2}
\end{array}
\end{vmatrix}
& \neq & 0\\
\end{eqnarray*}
We factor in each row but the last the corresponding $\vy$, resp. $\vi$
and add $d$ times the first column, and all following ones to the last:
\begin{eqnarray*}
\begin{vmatrix}
2(1+\pt x_n^{2d}) & 0 & \cdots & 2d\\
0 & (3d-1)(1+\pt x_n^{3d-1}) & \cdots & (3d-1)\\
\vdots & \vdots & \ddots & \vdots\\
-2\pt \vy y x_n^{2d-1} & -(3d-1)\pt \vi x_i x_n^{3d-2} & \cdots &
(3d-1)\vn+\vn' x_n^{3d-2}
\end{vmatrix}
& \neq & 0\\
\end{eqnarray*}
By hypothesis
each summand of the determinant is real and non-negative real,
and at least one summand is positive.
\qed


\begin{prop}
\labell{largest}
Suppose $\vn'\in\RR^{\geq0},\vy^3,\vi^\dd,\vn\in\rplus$.
Then for $\pt\in[0,\ptcrit[$
there is a unique critical point with $\vy y, \vi x_i\in\rplus,
x_n\in\rplus$, and it has maximal critical value.
\end{prop}

\proof
By assumption the coordinate change $\check y=\vy y$,
$\check x_i=\vi x_i$ transforms equations (\ref{crit-eqn})
into an equivalent system of equations:
\begin{eqnarray}
\check y^2 & = & \vy^3(1+\pt x_n^{2d})\notag
\label{eqn}
\\
\check x_i^{3d-1} & = & \vi^{3d}(1+\pt x_n^{3d-1})\\
x_n^{3d-1} & = & \vn + \vn' x_n^{3d-2}
+ \pt(2\check y x_n^{2d-1}
+(3d-1)\sum_i\check x_i x_n^{3d-2})\notag
\end{eqnarray}
For $\pt=0$ there is a unique solution with
$\check y,\check x_i, x_n\in \rplus$, because in that case
the last equation in (\ref{eqn}) has a unique
solution in $\rplus$ by the sign rule.

Since solutions depend continuously on the parameter $\pt$,
there are only the following transitions, which lead to a change of the
number of positive real solutions:
\begin{enumerate}
\item
a solution tends to infinity, which actually happens for the critical
parameter, but not before, cf.\ lemma \ref{simple-asym}
\item
positive real solutions become complex and vice versa,
\item
positive real solutions become semi-positive or vice versa, but
with the give hypotheses there is never a semipositive real solution,
since $x_n=0$ implies $\pt_n=0$ and
$\check y=0$ or $\check x_i=0$ implies $x_n=0$ by (\ref{eqn}).
\end{enumerate}
The uniqueness claim follows since also the second case can be
excluded:

Suppose a complex solution tends to a positive real solution, then
there is another complex solution tending to the same real solution,
since equations (\ref{eqn}) are real, thus complex conjugation acts on
solutions.
Therefore the limit solution is positive real and degenerate.
But with the given hypotheses such degenerate solutions do not exit
due to lemma \ref{hessian}.

To prove the maximality claim for the value of the distinguished
solution we notice:
\begin{enumerate}
\item
a critical point with last coordinate $x_n$ of smaller modulus than
the distinguished critical point has by (\ref{eqn}) all coordinates
smaller in modulus and hence smaller critical value by lemma
\ref{crit-val}.
\item
for $\pt=0$ the distinguished critical point
has  maximal last coordinate and uniquely so if $\vn'\neq0$.
\item
For $\pt_n'\in\rplus$ there is no second critical point with last
coordinate equal in modulus to that of the distinguished critical
point, since we may argue as follows: We get the last equation of 
(\ref{eqn}) for the last coordinate and for its absolute value,
\begin{eqnarray*}
x_n^{3d-1} & = & \vn + \vn' x_n^{3d-2}
+ \pt(2\check y x_n^{2d-1}
+(3d-1)\sum_i\check x_i x_n^{3d-2})\\
|x_n^{3d-1}| & = & \vn + \vn' |x_n^{3d-2}|
+ \pt(2\check y' |x_n^{2d-1}|
+(3d-1)\sum_i\check x'_i |x_n^{3d-2}|)
\end{eqnarray*}
where $\check y',\check x_i'$ denote the coordinates of the
distinguished critical point.

By the other equations in (\ref{eqn}) the absolute values of $\check y$
and $\check x_i$ are bounded by $\check y',\check x_i'$, hence
we may
deduce that all summand in the first equation are positive real.
In particular $x_n^{3d-1},x_n^{3d-2}\in\rplus$ and we conclude
$x_n\in\rplus$ and consequently $\check y,\check x_i\in\rplus$.
\end{enumerate}
By the continuity of critical points we may conclude, that for all
admissible $\pt$ and $\vn'$ the last coordinate of a critical
point is bounded by that of the distinguished critical point.
Therefore the distinguished critical points has maximal value
according to the first observation.
\qed

\subsection{asymptotic arcs, paths, and induced paths}

We recall that the parameter space $\CC^{N-1}$
of pencils of polynomials is naturally identified -- setting $x_0=0$ --
with a parameter space of polynomials contained
in $\CC[y,x_1,...,x_n]$ with vanishing constant coefficient.
Hence arcs and paths in $\CC^{N-1}$ are given by families of polynomials parameterised by a real
interval.

We consider arcs $\a_\jbold$ in the family $\gfami$
in bijection to tuples $j_1,...,j_{n-1}$
starting at the Hefez-Lazzeri base point
$\l_\ncount=v_\ncount,\,\pt_n'=0,\,\l=0,$
which corresponds to a polynomial $f=y^3-3 v_0 y+
\sum_\ncount^{n}
(x_\ncount^\dd-\dd v_\ncount x_\ncount)$ admitting a Hefez-Lazzeri
geometric basis.
They are composed of two parts each
\begin{enumerate}
\item
$\l,\pt_n'=0$, $\l_\ncount$ moves from $v_\ncount$
to $v_\ncount\xi^{j_\ncount-1}$ and $\pt_0$ from $v_0$ to
$v_0\xi^{d (j_0-1)}$,
cf. the construction of the paths $\oo_j$
in \reff{pathsandhomotopies}.
\item
$\pt=0,\pt_\ncount$ stay fixed, $\pt_n'$ increases to some 
small finite value.
\item
$\l_\ncount,\pt_n'$ stay fixed, $\l$ increases from $0$ to $\lcrit$.
\end{enumerate}
By lemma \ref{simple-asym}
each $\a_\jbold$ leads to a point of $\afami$
without intersecting $\afami$ elsewhere,
so we get well defined geometric elements associated to $\afami$.
An arbitrarily small isotopy yields the same geometric element
represented by a path $\cg_\jbold$ in the complement of
$\afami\cup\bfami$.

\begin{prop}
\labell{transport}
The relation imposed on the generators along the path $\cg_\jbold$
is 
$$
(t_{j1}\inv\d_0)^{\dd-1}=(\d_0t_{j1}\inv)^{\dd-1}.
$$
\end{prop}

\proof
From lemma \ref{leadc} we deduce that over a transversal disc
to $\afami$ the discriminant is an asymptotic germ of pole order
$\dd-1$.

So from prop.\ \reff{local1} we get a relation between elements
of a geometric basis in a local reference fibre
adapted to the asymptote.
The point is, that we have to identify the significant elements $b_0,b$
with elements in terms of the Hefez-Lazzeri geometric basis along
the arc $\a_\jbold$.

Now the difficult construction of this section shows, that a
representative of $t_j$ is first transported to a representative of
$t_{1..1}({\jbold})$ according to lemma \reff{hotopy}.
That representative is further transported to $b_0$, since its enclosed
puncture remains the puncture of largest modulus
according to proposition \ref{largest}
over the second and third part of $\a_\jbold$.

Of course a loop around all punctures is transported to
a loop around all punctures,
hence parallel transport identifies $\d_0$ with $b_0b$.
We conclude $b=t_j\inv\d_0$ and get
\begin{eqnarray*}
(t_j\inv\d_0)^{\dd-1}t_j & = & t_j(t_j\inv\d_0)^{\dd-1}\\
\iff
(t_j\inv\d_0)^{\dd-1} & = & \d_0(t_j\inv\d_0)^{\dd-2}t_j\inv,
\end{eqnarray*}
which yields the claim.
\qed

\begin{prop}
\labell{generators1}
The set of paths $\cg_\jbold$ generates $\pi_1(\CC^{N-1}-\afami)$.
\end{prop}

\proof
Since $\afami$ is defined by
$\ell_{n,d}=p_{n-1,d}^{\dd-1}\in\CC[u_\nu']=\CC[u_\nu|\nu_0=0]$,
the complement of $\afami$ projects to the complement of
$\dfami_{n-1,d}$ via linear projection along all $u_\nu,\nu_0>0$.
It thus suffices to prove that the projected paths
$\cg_\jbold'$
generate $\pi_1(\CC^{N'}-\dfami_{n-1,d})$.

We observe that under projection the first two pieces
of each path contract to a point,
but the crucial observation is, that the tail of each
$\cg'(\jbold)$
is homotopic to a path $t_{1\cdots1}(\jbold)$ conjugated by a path
from the base point of $\cg'(\jbold)$
to that of $t_{1\cdots1}(\jbold)$ given by
\begin{eqnarray*}
\varpi(\jbold):\quad
\l & \mapsto &
\left(\l v_0 \xi^{d(j_0-1)},
\l v_1\xi^{j_1-1},\cdots, \l v_{n-1}\xi^{j_{n-1}-1}, 1\right),
\end{eqnarray*}
where the entries correspond to the coefficients of the
monomials $x_\ncount$ and the constant term in the perturbation
of the Brieskorn Pham polynomial
$f'=y^3+\sum_{\ncount=1}^{n-1} x_\ncount^\dd$.

The claim then follows from two additional observations:
First -- as is seen by methods similar to those in
\ref{pathsandhomotopies} --
the paths $\oo_\jbold^*(t_{1\cdots1}(\jbold))$
generate $\pi_1(L_{v'}-\dfami'_L)$, where $v'=(v_0,v_1,...,v_{n-1})$
and $\dfami'=\dfami_{n-1}$.

Second, each $t_{1\cdots1}(\jbold)$ conjugated by $\varpi(\jbold)$
is homotopic to $\oo_\jbold^*(t_{1\cdots1}(\jbold))$ conjugated
by $\varpi(1\cdots1)$, 
since each $\varpi(\jbold)$ is homotopic to the concatenation of
$\oo_\jbold$ -- cf.\ \ref{pathsandhomotopies} --
and $\varpi(1\cdots1)$.
\qed

\section{moduli quotient}
\label{projective}

It is time now to recall that our interest is in the fundamental group of 
the moduli stack $\mnd$ of Weierstrass fibrations.
Of course we have first to define the appropriate moduli problem and
then to show that we can construct $\mnd$ as a quotient of $\und'$
by a suitable group action.

\paragraph{Definition}
Suppose $W_1, W_2$ are Weierstrass fibrations.
An isomorphism $\phi:W_1\to W_2$ is called an 
\emph{isomorphism of Weierstrass fibrations},
if $\phi$ preserves the section at infinity and fits into a commutative
diagram
\begin{eqnarray*}
W_1 & \stackrel\phi\tto & W_2\\
\downarrow & & \downarrow\\
\PP^n & \tto & \PP^n
\end{eqnarray*}

Over the parameter space $\und'$ there is the tautological Weierstrass
fibration $\wnd$. The factor group
$G=\CC^*\times\gl_{n+1}\big/_{\textstyle \CC^*}$ by the central
subgroup given by elements $\left(\l^{-d/2},\l1\!\!1\right)$
acts faithfully on $\wnd$ by the action induced from
$$
(\l,A)\cdot (\mathbf x, y_0,y,y_2)
\quad = \quad
(A\cdot\mathbf x, \l^2 y_0,y,\l\inv y_2).
$$

\begin{prop}
Every isomorphism of Weierstrass fibrations is induced from an
automorphism in $G$ of the tautological Weierstrass fibration $\wnd$.
\end{prop}

\proof
Given an isomorphism $\phi$ of smooth Weierstrass fibrations
we get immediately an induced automorphism of $\PP^n$
and we may pick some $A\in\gl_{n+1}$ inducing it.

Therefore we may assume without loss of generality, that $\phi$
induces the identiy on the base. Hence $\phi$ induces abstract
isomorphisms of plane cubic curves in Weierstrass normal form, 
mapping points at infinity to points at infinity.

Now it is well known, that each such isomorphism is induced by an
automorphism of $\PP^2$ of the form determined by
a non-vanishing complex number $\l$ acting on the coordinates
by $\l\cdot(y_0,y,y_2)=(\l^{2}y_0,y,\l\inv y_2)$.

Hence we get a map from $\PP^n$ to a subgroup of $\pgl_2$
isomorphic to $\CC^*$. Thus it must be constant and our claim is
proved.
\qed

According to this result it is natural to conceive the following definition:

\paragraph{Definition}
The moduli stack $\mnd$
of smooth Weierstrass fibrations is the quotient of the base space
$\und'$ by the induced action of $G$.

\begin{prop}
\labell{modquot}
There is an exact sequence
$$
\pi_1(\CC^*\times\gl_{n+1},1)
\tto\pi_1(\und, u),\tto\pi_1(\mnd,[u])\tto 1.
$$
where the action of $\CC^*$ is defined by $\l\cdot(p,q)=(\l^4p,\l^6q)$
and $\gl_{n+1}$ acts by linear coordinate change
$A\cdot(p,q)=(p\circ A,q\circ A)$.
\end{prop}

\proof
First we note that the affine group $\CC[x_0,...,x_n]_d$
acts on $\wnd$ by
$$
s(x)\cdot (\mathbf x,y_0,y,y_2)
\quad = \quad
(\mathbf x, y_0,y-s(x) y_0,y_2)
$$
and the induced action on $\und$ is free and faithful
and has $\und'$ as a transversal section. Hence
we may replace $\pi_1(\und,u)$ by $\pi_1(\und',u)$.
But then the claim is obviously just the exact sequence for
orbifold fundamental groups and the action of $G$ on $\und$
inducing the action of $G$ on $\und'$ made explicit.
\qed

In the following discussion
we are going to combine results on fundamental groups of $\und$ and various of
its subspaces which do \emph{not} have the same base point.

Still all these base points are contained in a ball in $\und$, so our convention is
that all occurring fundamental groups are identified using a connecting path
for their base points inside this ball, which makes the identification 
unambiguous.
\\[4mm]
Let $u$ now be the parameter point
corresponding to the Brieskorn-Pham hypersurface 
$$
u:\quad y^3+x_0^\dd+x_1^\dd+\cdots+x_n^\dd.
$$
Its $\CC^*$ orbit belong entirely to the affine Brieskorn-Pham
family $\ffami$ to which we shift our attention for the moment
$$
y^3+a_0x_0^\dd+a_1x_1^\dd+\cdots+a_nx_n^\dd.
$$

\paragraph{Definition}
The element $\delta_k$ is defined as the element represented by the path in the Brieskorn-Pham family
$y^3+\sum a_ix_i^\dd$ given by $a_i=1$ for $i\neq k$ and $a_k=e^{it}$.
\\

According to our remark above $\d_0$ may be identified with elements
represented by a loop in the coefficient $a_0=z$ of $x_0^\dd$
for any sufficiently small perturbation of the Brieskorn-Pham 
polynomial.

In particular $\d_0$ identifies with the element of the same name
from lemma \reff{cox}, and therefore can be expressed in the
geometric basis $t_i$ of $\pi_1(L_0-\dfami_L)$
using the enumeration function $\ibold_0:\{1,...,2(\dd-1)^n\}\to \ind$:
$$
\d_0\quad =\quad\prod_{k=1}^{2(\dd-1)^n} t_{\ibold_0(k)}.
$$
Remember that $\ibold_0$ has been defined above as the 
lexicographical order $\prec_0$ of 
$\ind$ derived from the \emph{reverse} order of the factors $\{1,2\}$,
$\{1,...,\dd-1\}$ of $\ind=\{1,2\}\times\{1,...,\dd-1\}^n$.
\begin{eqnarray*}
i_0 i_1\cdots i_n \prec_0 j_0 j_1 \cdots j_n
& \iff &
\exists k :\, i_\nu=j_\nu \forall \nu < k, i_k>j_k.
\end{eqnarray*}
To give analogous expressions for $\d_k$ we employ similar
enumeration functions
\begin{eqnarray*}
\ibold_\ncount &: & \{1,...,2(\dd-1)^n\}\to \ind
\end{eqnarray*}
which again are most conveniently described by the linear order
$\kord$ they induce on $\ind$:
\begin{eqnarray*}
i_0 i_1\cdots i_n \prec_\ncount j_0 j_1 \cdots j_n
& \iff &
 i_\ncount<j_\ncount\quad
 \vee\quad i_\ncount=j_\ncount,\:\:
i_0 i_1\cdots i_n \prec_0 j_0 j_1 \cdots j_n
\end{eqnarray*}
For $n=1$ in particular $\prec_1$ is the order on $\{1,2\}
\times\{1,...,d-1\}$ given by
$$
(2,1)\prec_1(1,1)\prec_1(2,2)\prec_1(2,1)\prec_1
\cdots\prec_1(2,\dd-1)
\prec_1(1,\dd-1).
$$

Let us resume now our argument to get expressions for
the $\d_k$:

\begin{lemma}
\labell{deltaab}
The $\delta_k$ all commute with each other and
$$
\prod\d_k^6\quad=\quad 1\in \pi_1(\und/\CC^*).
$$
\end{lemma}

\proof
The family $y^3+\sum a_ix_i^\dd$ has discriminant given by the 
normal crossing
divisor $\prod a_i$ and hence the fundamental group of the 
complement is abelian.

Since the $\d_k$ are geometric generators associated to the $n+1$ hyperplanes
the fundamental group of the quotient of the complement $\und$
by $\CC^*$ acting with multiplicity $6$ is the free abelian group
generated by the $\d_k$ modulo
the subgroup generated by the sixth power of their product.

Of course this relation maps homomorphically to $\pi_1(\und/\CC^*)$.
\qed

\begin{lemma}
\labell{ind-anf}
Consider the Hefez Lazzeri family
$$
y^3-v_0yx_0^{2d}
+x_1^\dd-v_1x_1x_0^{\dd-1}
+zx_0^\dd.
$$
Suppose $t_i$ and $t_i'$ form geometric
Hefez-Lazzeri bases for positive real $v_0> v_1$ respectively
$v_0'= v_1,v_1'=v_0$ of sufficiently distinct magnitude,
then there is a path connecting the base points such that the associated isomorphism
on fundamental groups is given by
$$
t_{i_0i_1}\quad \mapsto \quad t_{i_1i_0}'\quad 
$$
\end{lemma}

\proof
We first convince ourselves that $t_{11}=t'_{11}$ which follows immediately if we change $v_0,v_1$ continuously
in the real line swapping places since the extremal real puncture will keep that property and hence the
corresponding geometric element will not be changed.

To move $t_{i_0i_1}$ we first proceed along a path $\oo_{i_0i_1}$
as in lemma \ref{hotopy} so that it becomes the $t_{11}$ in the new 
system, then do the same as above and finally employ the
same path $\oo_{i_0i_1}$
back again to come to the final position.

The paths thus needed all connect the two base points with 
$v_0,v_1$ interchanged albeit in different ways.
But all possible concatenated paths are trivial in the complement
of the cuspidal bifurcation component.
They may still be non-trivial in the complement of the Maxwell
bifurcation component, but the induced isomorphism
on the fundamental group is trivial, since the corresponding braid
transformation just imposes the commutation relation, which is
needed to establish the isomorphism. 
\qed

\begin{lemma}
\labell{ind-fort}
Consider the Hefez Lazzeri family
$$
y^3-v_0yx_0^{2d}+
\sum (x_\ncount^\dd-v_\ncount x_\ncount x_0^{\dd-1})+zx_0^\dd.
$$
Suppose $\pi$ is a permutation such that $t_i$ and $t_i'$ are geometric generators
for positive real $v_0\gg v_1\gg\cdots\gg v_n$ respectively
$v'_{\pi(0)}=v_0,v'_{\pi(\ncount)}=v_\ncount$
of sufficiently distinct magnitude
then there is a path connecting the respective base points such that 
the associated isomorphism
on fundamental groups is given by
$$
t_{i_0i_1i_2\cdots i_n}\quad \mapsto \quad 
t_{i_{\pi(0)}i_{\pi(1)}i_{\pi(2)}\cdots i_{\pi(n)}}'
$$
\end{lemma}

\proof
The same idea of proof as above applies here.
\qed

\begin{lemma}
\labell{fermat}
The element $\d_1^\ffami$ in the Fermat family
$\ffami:\,a_0x_0^\dd+
a_1x_1^\dd$ can be expressed in the geometric basis $t_i$
as
$$
\d_1\quad=\quad t_1t_2\cdots t_{\dd-1}.
$$
\end{lemma}

\proof
Let $\ufami_{\PP^1,l}$ as in theorem \ref{zariski} denote the
discriminant complement associated to the complete linear
system of degree $l$ on $\PP^1$.
It is the quotient of the subset $\tilde\ufami_l\in\CC[x_0,x_1]_l$
of homogeneous polynomials of degree $l$ defining $l$
distinct points in $\PP^1$ modulo the diagonal $\CC^*$ action.
Hence there is an exact sequence
$$
1\tto \pi_1(\CC^*)\tto \pi_1(\tilde\ufami_l)\tto \pi_1(\ufami_{\PP^1,l})
\tto 1
$$
Now $\d_1\d_0\in \pi_1(\tilde\ufami)$ is the image of a generator
of $\pi_1(\CC^*)$ and 
with $\d_0\,=\, t_{l-1}\cdots t_1$ and 
Zariski's result (cf.\ theorem \ref{zariski})
$$
t_1\cdots t_{l-1}t_{l-1}\cdots t_1\quad=\quad 1
$$
we conclude $\d_1\,=\, t_1\cdots t_{\dd-1}$
for $l=\dd$ as claimed.
\qed

\begin{lemma}
\labell{anfang}
In case $n=1$ the element $\d_1$ can be expressed in the geometric basis $t_i$:
$$
\d_1\quad=\quad t_{2,1}t_{1,1}
t_{2,2}t_{1,2}\cdots t_{2,\dd-1}t_{1,\dd-1}.
$$
\end{lemma}

We postpone the proof to be able to use some ideas of the proof
of the following corollary:

\begin{lemma}
\labell{dim-ind}
In case of general $n$ the element $\d_1$ can be expressed in the geometric basis $t_i$
using the enumeration $\ibold_{1}:\{1,...,2(\dd-1)^n\}\to \ind$ as
$$
\d_1\quad=\quad \prod_{k=1}^{2(\dd-1)^n} t_{\ibold_{1}(k)}.
$$
\end{lemma}

\proof
Since the case $n=1$ is given in 
lemma \ref{anfang} we may assume inductively that the claim is
proved for $n-1$:
$$
\d_1^{(n-1)}\quad=\quad \prod_{k=1}^{2(\dd-1)^{n-1}} t_{i'_{1}(k)}.
$$
Hence there is a homotopy $H$ in the space
$\ufami_{n-1,d}$ between the representing path of $\d_1^{(n-1)}$
and the concatenation of the representing paths of
the $t_{i'_{1}(k)}$.

We consider the map defined on $\ufami_{n-1,d}$ by
$$
f'\quad \mapsto\quad f=f'+x_n^\dd
$$
which maps $H$ to $\und$.
Its image provides a homotopy between the representing
path of $\d_1^{(n)}$ and the concatenation of the images of the paths
representing the $t_{i'_{1}(k)}$.

These images represent elements $t^+_{i'_{1}(k)}$
which by lemma \ref{bundle} can be expressed as
$t^+_{i'_{1}(k)}=t_{i'_{1}(k)\,(d-1)}\cdots t_{i'_{1}(k)\,2}t_{i'_{1}(k)\,1}$.
Therefore the image of $H$ establishes the relation
$$
\d_1^{(n)} \quad = 
\prod_{k=1}^{2(\dd-1)^{n-1}} t^+_{i'_{1}(k)}\quad = 
\prod_{k=1}^{2(\dd-1)^{n-1}}
t_{i'_{1}(k)\,(d-1)}\cdots t_{i'_{1}(k)\,2}t_{i'_{1}(k)\,1}\quad = 
\prod_{k=1}^{2(\dd-1)^{n}} t_{\ibold_1(k)}.\\[-6mm]
$$
\qed

\proof[ of lemma \ref{anfang}]
Here we make use of a homotopy $H$ in $\ufami_{\PP^1,3d}$
between the representing paths of $\d_1^\ffami$ and the
concatenated path $t_1\cdots t_{\dd-1}$. 
We map the homotopy to $\ufami_{2,d}$ using the map
$f'\mapsto f=y^3+f'$.
By an argument as above we thus get
$$
\d_1 \quad = \quad
t_1^+\cdots t_{\dd-1}^+
\quad=\quad
t'_{1,2}t'_{1,1}\cdots t'_{\dd-1,2} t'_{\dd-1,1},
$$
where the $t'$ form a Hefez-Lazzeri base for a perturbation
with $v_0\ll v_1$. We therefore apply the lemma \ref{ind-anf} to
get
$$
\d_1 \quad = \quad
t_{2,1}t_{1,1}\cdots t_{2,\dd-1} t_{1,\dd-1},
$$
\qed

\begin{lemma}
\labell{proj-I}
In case of general $n$ an expression for $\d_\ncount$ is given by
$$
\d_\ncount
\quad=\quad \prod_{m=1}^{2(\dd-1)^n}t_{\ibold_\ncount(m)}.
$$
\end{lemma}

\proof
We get the expression for general $\d_\ncount$ using a transposition $\pi=(1\ncount)$ in lemma \ref{ind-fort} on the expression for $\d_1$.
Then it is obvious that the non-reversed order is transferred from the first entry into the $\ncount$-th entry.
\qed

\begin{prop}
\labell{projrel}
The image of a generator of $\pi_1(\CC^*)$ in the fundamental group $\pi_1(\und)$
for the natural map to the orbit of the base point is given in the Hefez-Lazzeri geometric generators
$t_i$ as
$$
\prod_{\ncount=0}^n\left(\prod_{m=1}^{2(\dd-1)^n}t_{\ibold_\ncount(m)}
\right)^6
$$
\end{prop}

\proof
Immediate from the previous.
\qed

\begin{lemma}
\labell{projact}
In the fundamental group $\pi_1(\und/\gl_n)$ the following relation
holds:
$$
\d_\ncount^\dd \quad = \quad 1.
$$
\end{lemma}

\proof
Each element $\d_\ncount$ is the trace of the Brieskorn-Pham
point transported by the $\CC^*$ action on the coefficient $a_\ncount$.
The $\CC^*$ action on the variable $x_\ncount$
has the same effect as the $\CC^*$ action with multiplicity $\dd$ on 
the coefficient of $x_\ncount^\dd$,
hence the claim.
\qed

\section{Conclusion}
\label{results}

Finally we are in the position to prove the new theorems of the
introduction with an actually more explicit claim:

Recall the definition of the index set
$\ind=\{\,(i_0,...,i_n)\,|\,1\leq i_\nu\leq \dd-1,i_0\leq2\,\}$,
and its reverse lexicographical order $\prec_0$.
We define a graph $\Gamma_{n,d}$ with vertex set $\ind$ and
edge set
$$
E(\Gamma_{n,d})\quad=\quad
\big\{(\ibold,\jbold)\,\big|\,\ibold\neq\jbold,\, i_\nu-j_\nu\in\{0,1\}
\,\forall\, \nu\,\big\}.
$$
Example of such graphs we have given in the introduction. The Main Theorems are thus
made precise with the definition of $\Gamma_{n,d}$ just given and the enumeration
functions $\ibold_\ncount$ of the previous section defining the distinguished elements
$\d_\ncount$

\paragraph{Main Theorem for discriminant complement}
\unitlength=1pt
The complement $\und$ of the discriminant $\dfami_{n,d}$
has fundamental group $\pi_1$ finitely presented by
generators $t_\ibold$, $\ibold\in \ind$ and relations
\begin{enumerate}
\item
$ t_\ibold t_\jbold= t_\jbold t_\ibold$,
for all $(\ibold,\jbold)\not\in E_{n,d}$,
\item
$ t_\ibold t_\jbold t_\ibold= t_\jbold t_\ibold t_\jbold$,
for all $(\ibold,\jbold)\in E_{n,d}$,
\item
$ t_\ibold t_\jbold t_\kbold t_\ibold= t_\jbold t_\kbold t_\ibold t_\jbold$,
for all $\ibold\prec\jbold\prec\kbold$
such that $(\ibold,\jbold),(\ibold,\kbold),(\jbold,\kbold)\in E_{n,d}$,
\item
for all $\ibold\in \ind$\\[-2mm]
\begin{equation}
\label{asym-rel}
 t_\ibold
 \left( t_\ibold\inv \prod_{\ibold\in\vernd^\prec} t_{\ibold}\right)^{\dd-1}
\,=\quad
\left( t_\ibold\inv \prod_{\ibold\in\vernd^\prec} t_{\ibold}\right)^{\dd-1}
t_\ibold
\end{equation}
\end{enumerate}

\paragraph{Main Theorem for moduli stacks}
\unitlength=1pt
The moduli stack $\mnd$
for Weierstrass hypersurfaces of type $(n,d)$
has fundamental group $\pi_1$ finitely presented by
generators $t_\ibold$, $\ibold\in \ind$ and relations
i), ii), iii) and iv) as above, and additionally v) and  vi):
\begin{enumerate}
\setcounter{enumi}{4}
\item
\begin{equation}
\label{Cact}
\delta_0^6\delta_1^6\cdots\delta_n^6 \quad=\quad1
\hspace*{8cm}
\end{equation}
\item
\begin{equation}
\label{Pact}
\delta_0^\dd  \quad=\quad1
\hspace*{9cm}
\end{equation}
\end{enumerate}

\proof[ of both theorems]
By lemma \ref{b-gen} there is a fibre $L'$ for some projection $p_v:\CC^{N-1}\to\CC^{N-2}$
such that lemma \ref{genab} applies;
elements $r_b$ of a geometric basis for $L'$ and elements $r_a$ generating $\pi_1(\CC^{N-1}-\afami)$
form a generating set for $\pi_1(\CC^{N-1}-\afami-\bfami)$.

So in an abstract sense a presentation is obtained
with the help of lemma \ref{pd},
but we still have to make explicit
the unspecified relations.

In fact we choose generators $t_\ibold,\ibold\in \ind$ represented by a
Hefez-Lazzeri geometric basis.
By prop.\ \ref{generators1} we may choose to take the generators
$r_a$ to be represented by the
paths $\cg_\jbold$. Then the relations of second type in
\ref{pd} can be replaced
using those of prop.\ \ref{transport}, so we get the relations in
(\ref{asym-rel}).

From prop.\ \ref{kernel} we infer that we may replace the first set of relations by any other
set normally generating the kernel of the map from the free group
to the fundamental group of the discriminant complement in the singularity unfolding,
in particular by the set given in below in the result cited
from \cite{habil}.

Proceeding now to the quotient of $\und'$ it suffices to refer to
proposition \ref{modquot} to
add just two relations, relation (\ref{Cact}), as we have shown
in prop.\ \ref{projrel} and relation (\ref{Pact}) due to lemma
\ref{projact}.
\qed

\begin{thm}
\labell{singularity}
The fundamental group of the discriminant complement for any versal
unfolding of a Brieskorn Pham polynomial
$y^3+x_1^\dd+\cdots x_n^\dd$ 
is finitely presented by
generators $t_\ibold$, $\ibold\in \ind$ and relations
\begin{enumerate}
\item
\hspace*{4.8mm}
$t_\ibold t_\jbold= t_\jbold t_\ibold$,
\hspace*{5.1mm}
for all $\ibold,\jbold$ such that $(\ibold,\jbold)\not\in E_{n,d}$,
\item
\hspace*{2mm}
$ t_\ibold t_\jbold t_\ibold= t_\jbold t_\ibold t_\jbold$,
\hspace*{2.1mm}
for all $\ibold,\jbold$ such that $(\ibold,\jbold)\in E_{n,d}$,
\item
$ t_\ibold t_\jbold t_\kbold t_\ibold= t_\jbold t_\kbold t_\ibold t_\jbold$,
for all $\ibold\prec\jbold\prec\kbold$
such that $(\ibold,\jbold),(\ibold,\kbold),(\jbold,\kbold)\in E_{n,d}$,
\end{enumerate}
\end{thm}

\subsection{special cases for small invariants}

\paragraph{The case of moduli stack of elliptic curves}
In case $n=0$ we consider the set $\ufami'$ of three distinct points
in $\CC$ with center of mass in $0$.

Accordingly the moduli stack is that of elliptic curves and our theorem
just reproduces the well-know result, that
\begin{eqnarray*}
\pi_1(\mfami_{0})
&\cong&
\langle\:t_1,t_2\:|\: t_1t_2t_1\,=\, t_2t_1t_2,\quad
\d_0^6=(t_1t_2)^6=1\:\rangle\\
&\cong&
\slz.
\end{eqnarray*}

\paragraph{The case of elliptic curves in the Hirzebruch surface
$\FF_1$}
In case $n=1,d=1$ we consider the set $\ufami_{1,1}$
of curves in $\FF_1$,
which is the subset of all smooth elliptic curves in the class of a triple
positive section.

That set can be understood as the set of elliptic curves in $\PP^2$ not
passing through a distinguished point (the blow down of the
exceptional line in $\FF_1$).
Therefore its $\pi_1$ should be a central extension of $\pi_1$
of the set of all smooth elliptic curves in $\PP_2$.

Our results yield a finite presentation of $\pi_1(\ufami_{1,1})$
by generators $t_{11},t_{12},t_{21},t_{22}$ and relations
\begin{itemize}
\item
$t_2t_3=t_3t_2$,
\item
$t_it_jt_i=t_jt_it_j$ if $(ij)\in\{(12),(13),(24),(34)\}$,
\item
$t_it_jt_kt_i=t_jt_kt_it_j$ if $(ijk)\in\{(124),(134)\}$,
\item
$t_4t_3t_2t_4t_3t_2t_1=t_1t_4t_3t_2t_4t_3t_2$,
\item
$t_3t_2t_1t_3t_2t_1t_4=t_4t_3t_2t_1t_3t_2t_1$,
\end{itemize}

In fact $\pi_1(\ufami_{1,1})$
is a central extension of
$\pi_1(\ufami_{\PP^2,3})$
which is finitely presented with the same generators and relations
except for one additional relation
\begin{itemize}
\item
$t_4t_3t_2t_1t_2t_1t_4t_3t_3t_1t_4t_2=1$
\end{itemize}

Moreover $\pi_1(\ufami_{1,1})$ is isomorphic to
$\pi_1(\tilde\ufami_{\PP^2,3})$, for the affine cone
$\tilde\ufami_{\PP^2,3}$ of $\ufami_{\PP^2,3}$
and the central extension is analogous to that
of the proof of lemma \ref{fermat}.

\paragraph{The case of smooth rational Weierstrass fibrations}
In case $n=1,d=2$ we consider the moduli stack $\mfami_{1,2}$
of smooth Weierstrass fibrations which are rational.
In fact every rational elliptic fibration with a section has a unique
Weierstrass model, which is smooth if and only if all fibres are
irreducible.
Hence our fundamental group is also that of the moduli stack
of rational elliptic fibrations with a section and irreducible fibres only.
$\pi_1(\mfami_{1,2})$ is finitely presented by generators
\begin{eqnarray*}
t_{1,1},t_{1,2},t_{1,3},t_{1,4},t_{1,5},t_{2,1},t_{2,2},t_{2,3},t_{2,4},
t_{2,5},
\end{eqnarray*} 
subjected to relations
\begin{eqnarray*}
&t_{i_0,i_1}t_{j_0,j_1}t_{i_0,i_1}
=t_{j_0,j_1}t_{i_0,i_1}t_{j_0,j_1}
&
\text{if } |i_1-j_1|\leq1,
(i_0-j_0)(i_1-j_1)\geq0\\
&t_{i_0,i_1}t_{j_0,j_1}
=t_{j_0,j_1}t_{i_0,i_1}
&
\text{if } |i_1-j_1|\geq2\text{ or }
(i_0-j_0)(i_1-j_1)<0\\
&\hspace*{3mm}t_{1,i_1}t_{2,i_1+1}t_{1,i_1+1}t_{1,i_1}
=t_{1,i_1+1}t_{1,i_1}t_{2,i_1+1}t_{1,i_1+1}\\
&t_{1,i_1}t_{2,i_1+1}t_{2,i_1}t_{1,i_1}
=t_{2,i_1}t_{1,i_1}t_{2,i_1+1}t_{2,i_1}\\
&(t_i\inv\d_0)^5 = (\d_0 t_i\inv)^5\\
&\d_0^6=1=\d_1^6
\end{eqnarray*}




\begin{thebibliography}{ACT}  





\bibitem[Be]{bessis} Bessis, David:
{\sl Zariski theorems and diagrams for braid groups},
Invent. Math. 145 (2001), 487--507






\bibitem[DL]{dl} Dolgachev, Igor, Libgober, Anatoly: {\sl On the fundamental group of the complement to
a discriminant variety}, in {\sl Algebraic geometry, Chicago 1980}, LNM 862,
Springer, Berlin-New York, 1981, 1--25,

\bibitem[D]{d} Dimca, Alexandru: {\sl Singularities and Topology of Hypersurfaces},
Springer-Verlag, New York, 1992, xvi+263 pp

\bibitem[FvB]{fvb} Fadell, Edward, Van Buskirk, James: {\sl On the braid groups of $E\sp{2}$ and $S\sp{2}$},
Bull. AMS 67 (1961), 211--213,

\bibitem[HL]{hl}
A.~Hefez, F.~Lazzeri:
{\sl The intersection matrix of Brieskorn singularities},
Invent.\ Math.\ 25 (1974), 143--157.




\bibitem[L\"o1]{reg} M. L\"onne: {\sl Monodromy groups of regular
elliptic surfaces},
Math.\ Z. 239  (2002),  no. 3, 441--453

\bibitem[L\"o2]{ireg} M. L\"onne: {\sl Monodromy groups of irregular
elliptic surfaces},
Compositio Math. 133 (2002), no. 1, 37--48

\bibitem[L\"o3]{habil} M. L\"onne: {\sl Braid monodromy of hypersurface singularities},
Habilitations\-schrift (2003), Hannover, mathAG/0602371

\bibitem[L\"o4]{bifbraid} M. L\"onne: {\sl Bifurcation braid monodromy 
of elliptic surfaces},
Topology, mathAG/0602371

\bibitem[L\"o5]{pdisc} M. L\"onne: {\sl Fundamental groups of spaces
of smooth hypersurfaces},
preprint, mathAG/0602371


\bibitem[Mi]{mi} Miranda, Rick: {\sl The moduli of Weierstrass fibrations
over $\PP^1$}, Math.\ Ann.\ {\bf 255} (1981), 379--394

\bibitem[Sei1]{sei1} Seiler, Wolfgang: {\sl Global moduli for elliptic 
surfaces with a section}. Compositio Math.\ {\bf 62} (1987), 169--185

\bibitem[Sei2]{sei2} Seiler, Wolfgang: {\sl Global moduli for polarized
elliptic surfaces}. Compositio Math.\ {\bf 62} (1987), 187--213

\bibitem[Za]{z} Zariski, Oscar: {\sl On the Poincar\'e group of rational plane curves},
Am. J. Math. 58 (1936), 607-619

\end{thebibliography}
\end{document}